\documentclass[12pt]{amsart}
\usepackage[margin=1in]{geometry}  
\usepackage{graphicx}              
\usepackage{amsmath}               
\usepackage{amsfonts}              
\usepackage{amsthm}                
\usepackage{caption} 
\usepackage{subcaption}
\usepackage{tikz}
\usepackage{enumitem}

\newtheorem{theorem}{Theorem}[section]
\newtheorem{lemma}[theorem]{Lemma}

\newtheorem{proposition}[theorem]{Proposition}
\newtheorem{corollary}[theorem]{Corollary}

\theoremstyle{remark}
\newtheorem{remark}[theorem]{Remark}

\theoremstyle{observation}

\theoremstyle{definition}
\newtheorem{definition}[theorem]{Definition}

\begin{document}

\title{Symplectically Replacing Plumbings with Euler Characteristic 2 4-Manifolds}
\author{Jonathan Simone}

\maketitle 

\begin{abstract}
We introduce new symplectic cut-and-paste operations that generalize the rational blowdown. In particular, we will define $k$-replaceable plumbings to be those that, heuristically, can be symplectically replaced by Euler characteristic $k$ 4-manifolds. We will then classify 2-replaceable linear plumbings, construct 2-replaceable plumbing trees, and use one such tree to construct a symplectic exotic $\mathbb{C}P^2\#6\overline{\mathbb{C}P^2}$.
\end{abstract}

\section{Introduction}

In recent years, symplectic cut-and-paste operations have been used to construct (symplectic) exotic 4-manifolds with ``small" $b_2$. To perform such an operation, one must:
\begin{itemize}
\item find a symplectic 4-manifold $(P,\omega_1)$ with strongly convex boundary embedded in an ambient symplectic 4-manifold $(X,\omega)$; 
\item construct a 4-manifold $B$ such that $\partial B=\partial P$ and such that $B$ admits a symplectic structure $\omega_2$ with strongly convex boundary; and
\item ensure that the induced contact structures on $\partial P$ and $\partial B$ are contactomorphic.
\end{itemize}

\noindent If these conditions are met, then by a result of Etnyre \cite{etnyreconvexity}, $Z=(X-\text{int}(P))\cup B$ inherits a symplectic structure from $\omega$ and $\omega_2$. Note that as a smooth 4-manifold, $Z$ may depend on the choice of contactomorphism. If the induced contact structures happen to be isotopic, however, then by choosing a contactomorphism isotopic to the identity, $Z$ is well-defined as a smooth 4-manifold.

\begin{definition} Suppose $P$ and $B$ are 4-manifolds with $\partial P=\partial B$. If $P$ and $B$ admit symplectic structures with strongly convex boundary that induce isotopic contact structures, and if $P$ is embedded in an ambient symplectic 4-manifold $X$, then we say that $P$ can be \textit{symplectically replaced} by $B$ and we call $B$ a \textit{symplectic replacement} of $P$.\end{definition}

An oft-used symplectic cut-and-paste operation is the rational blowdown, in which a negative-definite plumbing of $D^2$-bundles over $S^2$ is excised from a 4-manifold and a rational homology ball is glued in its place. This operation has been used to construct (symplectic) exotic 4-manifolds with small $b_2$. For example, in \cite{park}, Park constructed an exotic $\mathbb{C}P^2\#7\overline{\mathbb{C}P^2}$ and in \cite{stipsicz}, Stipsicz and Szab\'o constructed exotic $\mathbb{C}P^2\#6\overline{\mathbb{C}P^2}$s. 

The rational blowdown was introduced for linear plumbings by Fintushel and Stern \cite{fintushelstern}, generalized by Park \cite{parkbdown}, and shown to be symplectic by Symington \cite{symington2}. In \cite{sswahl}, Stipsicz-Szab{\' o}-Wahl generalized the operation to plumbing trees. Combining Park's definition of the rational blowdown with a result of Lisca \cite{lisca}, the linear plumbings that can be rationally blown down are precisely those that can be symplectically replaced by rational balls. Moreover, these are precisely the plumbings whose lens space boundaries $L(p,q)$ satisfy $p=m^2$ and $q=mn-1$, where $m>n>0$ are coprime integers. Note that these lens spaces were already known to (smoothly) bound rational balls by Casson and Harer \cite{cassonharer}. 

Equivalently, the plumbings that can be rationally blown down are those that can be obtained by the following inductive procedure: if the linear plumbing with framings $(-b_1,\ldots,-b_k)$  can be rationally blown down, then the plumbings with framings $(-b_1-1,\ldots,-b_k,-2)$ and $(-2,-b_1,\ldots,-b_k-1)$ can also be rationally blown down. The first such plumbing is the $-4$-disk bundle over $S^2$. This operation will arise many times throughout the paper, so we give it a name.

\begin{definition} Let $P$ be a linear plumbing with weights $(-b_1,\ldots,-b_k)$, where $b_i\ge2$ for all $i$. The \textit{buddings} of $P$ are the linear plumbings with weights $(-2,-b_1,\ldots,-b_{k-1},-b_k-1)$ and $(-b_1-1,-b_2,\ldots,-b_k,-2)$.\end{definition}

Since rational homology balls have Euler characteristic 1, a natural generalization of the rational blowdown is the following.

\begin{definition} A negative-definite plumbing $P$ is called \textit{$k$-replaceable} if it can be symplectically replaced by a negative-definite, minimal symplectic 4-manifold $B$ satisfying $\chi(B)=k$ and $b_3(B)=0$. We say that $P$ can be \textit{k-replaced} by $B$ and we call $B$ a \textit{k-replacement} of $P$.\end{definition}

Notice that 1-replaceable plumbings are precisely those that can be rationally blown down. Our goal is to use $k$-replaceable plumbings to construct closed, simply connected, symplectic, exotic 4-manifolds with small $b_2$. Thus we would like $B$ to be an Euler characteristic $k$ manifold with the smallest possible second Betti number. This is why we require that $B$ is minimal and that $b_3(B)=0$. We further require $B$ to be negative-definite so that Michalogiorgaki's gluing formula in \cite{micha} for Seiberg-Witten invariants is applicable. Moreover, by considering the long exact sequences of the pairs $(B,\partial B)$ and $(P,\partial P)$, since $B$ and $P$ are negative-definite and $b_3(P)=b_3(B)=0$, it follows that $b_1(P)=b_1(B)$.

In this paper, we will be mainly concerned with 2-replaceable plumbings. The first result is a classification of 2-replaceable linear plumbings. Note that there are infinitely many linear plumbings that have Euler characteristic 2, namely the disk bundles over $S^2$. We call such plumbings \textit{trivially} 2-replaceable. Starting with one of these disk bundles we can easily construct infinitely many 2-replaceable linear plumbings by plumbing these disk bundles with 1-replaceable linear plumbings (since the latter can be rationally blown down, reducing the Euler characteristic to 2).  These are the plumbings shown in Theorem \ref{thm:main}$(a)$. Since these are easy to construct, we are more interested in families of 2-replaceable linear plumbings that are not of this form, such as those in Theorem \ref{thm:main}$(b)-(d)$. 

\begin{theorem} Let $(-b_1,\ldots,-b_k)$ and $(-c_1,\ldots,-c_l)$ be obtained by sequences of buddings of $(-4)$ and let $z\ge2$ be any integer. Then a minimal linear plumbing is 2-replaceable if and only if it is either of the form:
\vspace{.1cm}
\begin{enumerate}[label=(\alph*)]
\item \includegraphics[scale=.6]{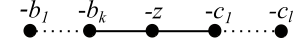} \qquad\qquad\qquad\qquad  for $k,l\ge0$
\end{enumerate}
\smallskip
or can be obtained by a sequence of buddings of one of the linear plumbings of the form:
\vspace{.1cm}
\begin{enumerate}[resume, label=(\alph*)]
\item \includegraphics[scale=.6]{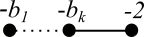}\quad (or \quad \includegraphics[scale=.6]{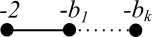}) \qquad\quad for $k\ge0$.
\vspace{.1cm}
\item \includegraphics[scale=.6]{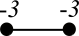}
\vspace{.1cm}
\item \includegraphics[scale=.6]{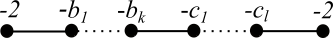} \qquad\qquad\quad for $k,l\ge1$
\end{enumerate}
\label{thm:main}\end{theorem}
\bigskip

\begin{remark} The proof of Theorem \ref{thm:main} relies on Lisca's classification (in \cite{lisca}) of symplectic fillings of lens spaces equipped with the canonical contact structure inherited from the unique tight contact structure on $S^3$. Thus, in principle, the theorem answers the question ``Which lens spaces (equipped with the standard contact structure) have strong symplectic fillings of Euler characteristic 2?" Moreover, the proof can be adapted to find (and classify) families of $k$-replaceable linear plumbings for $k\ge 3$. \end{remark}


\begin{theorem} For any integers $n,m\ge3$, the following are families of 2-replaceable trees:\\
\begin{center}\includegraphics[scale=.6]{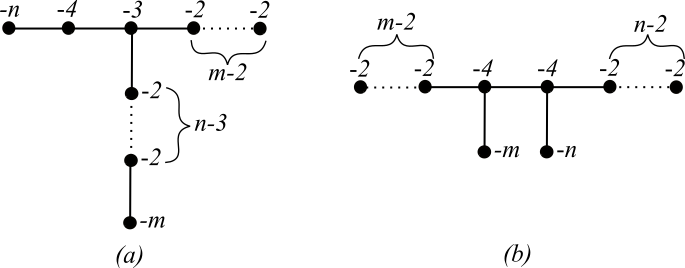}\end{center}
\label{thm:trees}
\end{theorem}

\begin{remark}\normalfont The families of plumbing trees in Theorem \ref{thm:trees} will be constructed from the (2-replaceable) linear plumbing with weights $(-2,-4,-4,-2).$ In the proof of the theorem, we will show that this linear plumbing is indeed 2-replaceable without relying on Theorem \ref{thm:main}. Instead, we will apply the theory of Lefschetz fibrations. It turns out that the families of plumbing trees of Theorem \ref{thm:trees} are interesting in the sense that they cannot all be built trivially by plumbing the 1-replaceable trees of \cite{sswahl} to a disk bundle over $S^2$ (c.f. the plumbings in Theorem \ref{thm:main}a). Moreover, the technique used in the proof of Theorem \ref{thm:trees} can be applied to obtain more families of 2-replaceable trees. For example, instead of starting with the linear plumbing with weights $(-2,-4,-4,-2)$, one could start with a different 2-replaceable linear plumbing.\end{remark}

Finally, using the 2-replaceable tree of Theorem \ref{thm:trees}(a) with $n=9$ and $m=3$ we perform symplectic cut-and-paste to construct the following.

\begin{theorem} The 2-replaceable tree of Theorem \ref{thm:trees}(a) with $n=9$ and $m=3$ can be embedded in $\mathbb{C}P^2\#16\overline{\mathbb{C}P^2}$. Call this tree $P$ and let $B$ denote its Euler characteristic 2 replacement. Then $X=(\mathbb{C}P^2\#16\overline{\mathbb{C}P^2}\setminus \textup{int}(P)) \cup_{\partial P}B$ is homeomorphic but not diffeomorphic to $\mathbb{C}P^2\#6\overline{\mathbb{C}P^2}$. Furthermore, $X$ admits a symplectic structure.\label{thm:exotic}\end{theorem}

This paper is organized as follows. In Section \ref{pf1} we will use Lefschetz fibrations and a lemma due to Endo-Mark-Van-Horn Morris \cite{mark} to prove Theorem \ref{thm:trees}. In Section \ref{exotic} we will use symplectic cut-and-paste to construct the exotic $\mathbb{C}P^2\#6\overline{\mathbb{C}P^2}$ of Theorem \ref{thm:exotic}. In Section \ref{contfracsection}, we will prove some facts about Hirzenbruch-Jung continued fractions that are needed for the proof of Theorem \ref{thm:main}, which can be found in Section \ref{mainproof}.

\subsection*{Acknowledgements} I am very grateful to my advisor Thomas E. Mark for providing me with the initial question that led to this work and for his help, support, and patience along the way.

\section{Proof of Theorem \ref{thm:trees}}\label{pf1}

\subsection{Lefschetz fibrations and the Key Lemma}\label{kl}

In this section, we will highlight the strategy used to prove Theorem \ref{thm:trees}. We assume the reader is familiar with Lefschetz fibrations and open book decompositions. Let $P$ be a symplectic negative-definite plumbing with strongly convex boundary that admits a symplectic Lefschetz fibration over $D^2$ with monodromy $\tau$ that can be written down in an explicit factorization. This monodromy naturally describes an open book decomposition of $Y$ that supports the contact structure $\xi$ induced by the symplectic structure. 
Suppose there is a different factorization of $\tau$ into right Dehn twists about homologically essential curves such that the total space $B$ of the corresponding Lefschetz fibration has Euler characteristic 2. Then, $B$ is a symplectic replacement of $P$. Since the obvious handlebody diagram of $B$ obtained from the monodromy has no 3-handles, we have that $b_3(B)=0$. Finally, in \cite{etnyre}, Etnyre showed that any strong symplectic filling of a contact manifold supported by a planar open book is negative-definite. Thus $B$ is a 2-replacement of $P$.

We will apply the following \textit{Key Lemma} due to Endo, Mark, and Van Horn-Morris in \cite{mark} to the monodromy factorizations associated to $P$ and $B$.

\begin{figure}[h!]
\centering
\includegraphics[scale=.4]{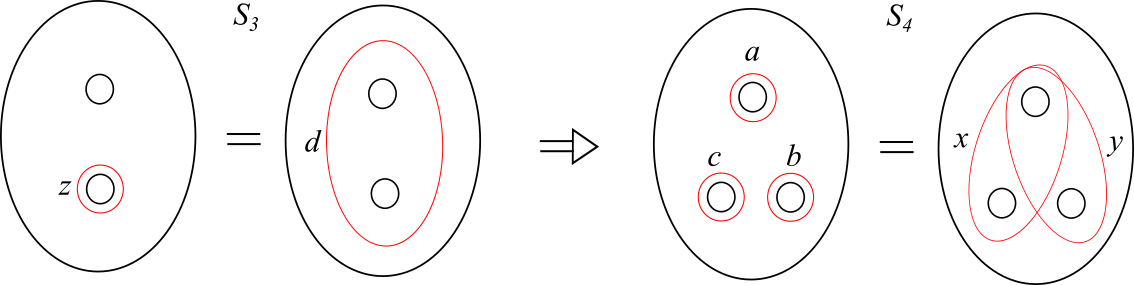}
\caption{The Key Lemma}\label{keylem}
\end{figure}

\begin{lemma} (Key Lemma \cite{mark}) Let $F$ be a planar surface containing as a subsurface a pair of pants, $S_3$. Let $z$ and $d$ be the boundary parallel curves marked in Figure \ref{keylem} and let the boundary component of $S_3$ corresponding to $z$ coincide with a component of $\partial F$. Let $F'$ be the planar surface obtained from $F$ by gluing a disk with two holes into the hole enclosed by $z$. Suppose that in the planar mapping class group $Mod(F,\partial F)$, the relation $w_1zw_2=w_1'dw_2'$ holds for some $w_1,w_2,w_1',w_2'\in Mod(F,\partial F)$. If $a$ commutes with either $w_1$ and $w_1'$ or $w_2$ and $w_2'$, then in $Mod(F',\partial F')$ we have the relation $w_1abcw_2=w_1'xyw_2'$.\label{keylemma}\end{lemma}

Assume the Key Lemma applies to the monodromies of $P$ and $B$ and suppose $P$ contains the curve $z$ and $B$ contains the curve $d$, as depicted in Figure \ref{keylem}. Let  $P'$ and $B'$ denote the total spaces the Lefschetz fibrations associated to the two new equivalent monodromy factorizations obtained from the Key Lemma. Then $Y'=\partial P'=\partial B'$ has an open book decomposition that can be described by these two factorizations. By Giroux's correspondence (\cite{girouxcorresp}), $Y=\partial P'=\partial B'$ admits a contact structure $\xi'$ that is supported by both open books. By \cite{gaymark} and \cite{plamenevskaya}, $P'$ and $B'$ both admit symplectic structures that are (strong) symplectic fillings of $(Y',\xi')$. Thus $B'$ is a symplectic replacement of $P'$.

We now claim that $\chi(B')=2$, $b_3(B')=0$, and that $B'$ is negative-definite. Since the obvious handlebody diagram of $B'$ obtained from its monodromy has one more 1-handle and one more 2-handle than the obvious handlebody diagram of $B$, it follows that $\chi(B')=\chi(B)=2$. Since there are no 3-handles in this diagram, we have that $b_3(B')=0$. Once again, by a result of Etnyre in \cite{etnyre}, since $B'$ is a strong symplectic filling of a contact manifold supported by a planar open book, $B'$ is negative-definite. 

Finally, for $B'$ to be a 2-replacement of $P'$, $B'$ must be minimal. We restrict our attention to the plumbings of Theorem \ref{thm:trees}. Using the above arguments, we will construct these plumbings and Euler characteristic 2 symplectic replacements in the next section. Suppose $P'$ is such a plumbing and let $B'$ be its Euler characteristic 2 symplectic replacement. If $B'$ is not minimal, then we can symplectically blow down a symplectic $-1$-sphere to obtain a 1-replacement of $P'$. In other words, $P'$ can be symplectically rationally blown down. All such plumbing trees are classified by Stipsicz-Szab{\'o}-Wahl in \cite{sswahl}. Since $P'$ is not among those trees, $B'$ must be minimal and so $P'$ is 2-replaceable.


\subsection{Proof of Theorem \ref{thm:trees}}\label{treesproof}

To construct the families of plumbing trees of Theorem \ref{thm:trees}, we will iteratively apply the Key Lemma. By the remarks above, these trees will automatically be 2-replaceable. All monodromy factorizations will be products of right Dehn twists around simple closed curves. For simplicity, a curve and a right Dehn twist about the curve will have the same label.

\begin{figure}[h!]
\centering
\includegraphics[scale=.4]{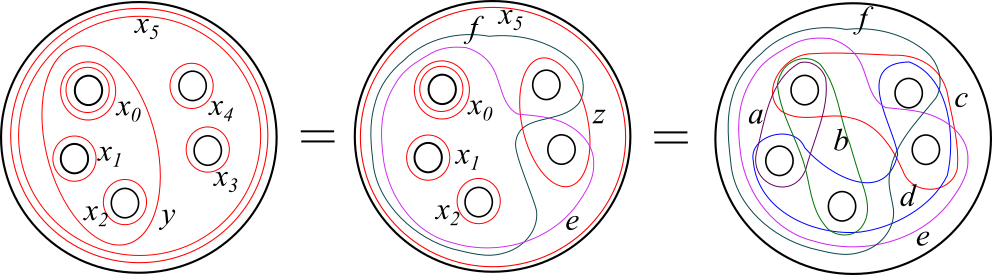}
\caption{$x_0^2x_1x_2x_3x_4x_5^2y=x_0^2x_1x_2zx_5ef=abcdef$}\label{relation}
\end{figure}

Let $P$ be the linear plumbing with framings $(-2,-4,-4,-2)$. $P$ can be viewed as a Lefschetz fibration with the monodromy factorization drawn on the left side of Figure \ref{relation}. It is given by $x_0^2x_1x_2x_3yx_4x_5^2$. Using a lantern relation applied to $x_3yx_4x_5$, we obtain the middle factorization in Figure \ref{relation}, $x_0^2x_1x_2zefx_5=x_0^2x_1x_2zx_5ef$. Finally, using the more general \textit{daisy relation} (defined in \cite{mark}), applied to $x_0^2x_1x_2zx_5$, we obtain the factorization $abcdef$ pictured on the right side of Figure \ref{relation}. By drawing a handlebody diagram of the total space of the Lefschetz fibration described by the monodromy factorization $abcdef$, easy homology calculations show that this 4-manifold is a 2-replacement of $P$ (see, for example, \cite{mark}).

\begin{figure}[h!]
\centering
\begin{subfigure}[h!]{\linewidth}
\centering
\includegraphics[scale=.4]{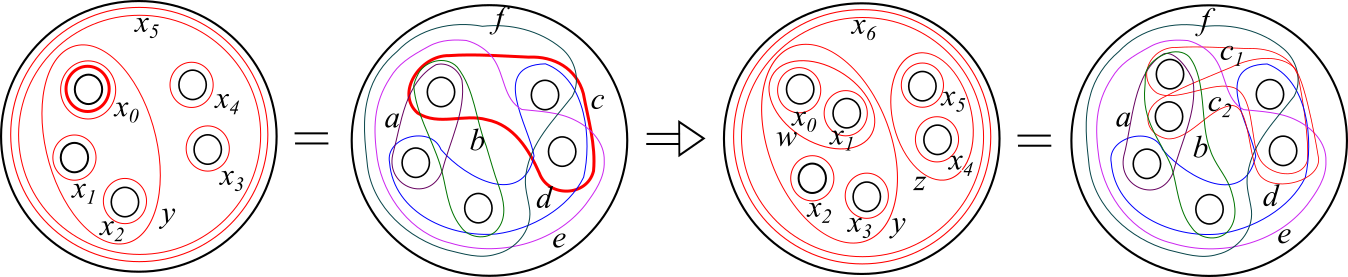}
\caption{$x_0^2x_1x_2x_3x_4x_5^2y=abcdef$ implies $x_0x_1x_2x_3x_4x_5x_6^2wyz=abc_1c_2def$}\label{keylemma1}
\end{subfigure}
\begin{subfigure}[h!]{\linewidth}
\centering
\includegraphics[scale=.4]{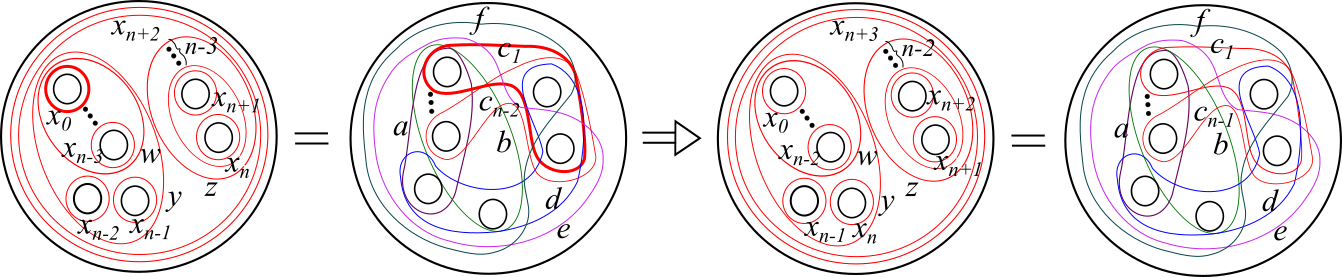}
\caption{$x_0\cdot\cdot\cdot x_{n+2}^2wyz^{n-3}=abc_1\cdot\cdot\cdot c_{n-2}def$ implies $x_0\cdot\cdot\cdot x_{n+3}^2wyz^{n-2}=abc_1\cdot\cdot\cdot c_{n-1}def$}\label{keylemma2} 
\end{subfigure}
\begin{subfigure}[h!]{\linewidth}
\centering
\includegraphics[scale=.4]{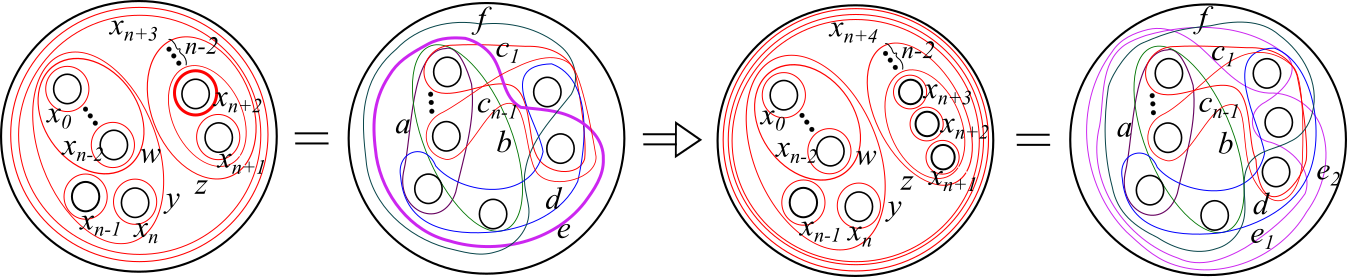}
\caption{$x_0\cdot\cdot\cdot x_{n+3}^2wyz^{n-2}=abc_1\cdot\cdot\cdot c_{n-1}def$ implies $x_0\cdot\cdot\cdot x_{n+3}x_{n+4}^3wyz^{n-2}=abc_1\cdot\cdot\cdot c_{n-1}de_2e_1f$}\label{keylemma3}
\end{subfigure}
\begin{subfigure}[h!]{\linewidth}
\centering
\captionsetup{justification=centering}
\includegraphics[scale=.4]{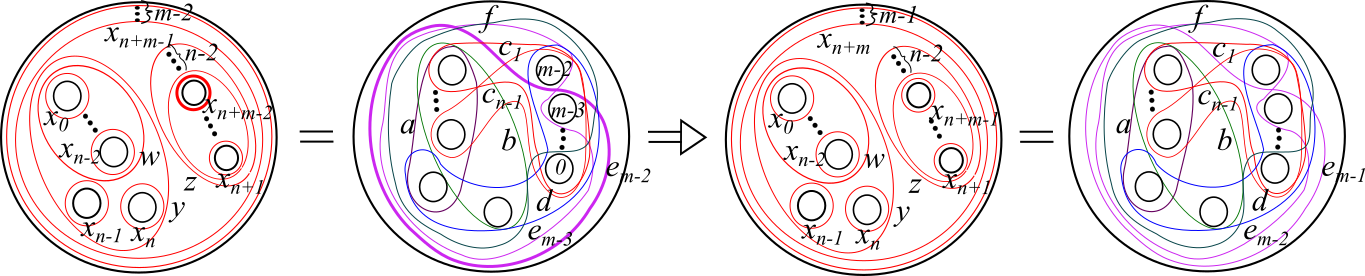}
\caption{$x_0\cdot\cdot\cdot x_{n+m-2}x_{n+m-1}^{m-2}wyz^{n-2}=abc_1\cdot\cdot\cdot c_{n-1}de_{m-2}\cdot\cdot\cdot e_1f$ implies\\ $x_0\cdot\cdot\cdot x_{n+m-1}x_{n+m}^{m-1}wyz^{n-2}=abc_1\cdot\cdot\cdot c_{n-1}de_{m-1}\cdot\cdot\cdot e_1f$}\label{keylemma4} 
\end{subfigure}
\caption{Repeated applications of the Key Lemma to the bold circles}\label{keylemma1-4}
\end{figure}

\begin{figure}[h!]
\centering
\begin{subfigure}[h!]{.4\textwidth}
\centering
\includegraphics[scale=.6]{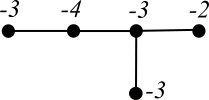}
\caption{Plumbing associated to Figure \ref{keylemma1}}\label{plumb1}
\end{subfigure}
\begin{subfigure}[h!]{.5\textwidth}
\centering
\includegraphics[scale=.6]{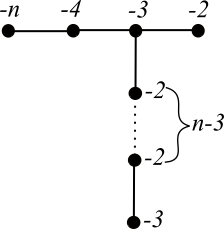}
\caption{Plumbing associated to Figure \ref{keylemma2}}\label{plumb2}
\end{subfigure}
\begin{subfigure}[h!]{.4\textwidth}
\centering
\includegraphics[scale=.6]{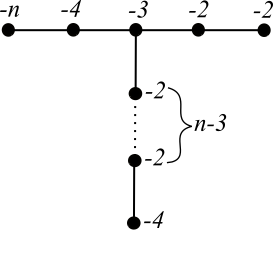}
\caption{Plumbing associated to Figure \ref{keylemma3}}\label{plumb3}
\end{subfigure}
\begin{subfigure}[h!]{.5\textwidth}
\centering
\includegraphics[scale=.6]{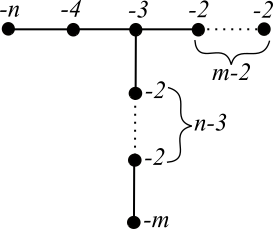}
\caption{Plumbing associated to Figure \ref{keylemma4}}\label{plumb4}
\end{subfigure}
\caption{2-replaceable plumbings associated to the monodromies in Figure \ref{keylemma1-4}}\label{plumb1-4}
\end{figure}

Now we will repeatedly apply the Key Lemma to the relation $x_0^2x_1x_2x_3x_4x_5^2y=abcdef$ to build a family of 2-replaceable trees. Notice that  the curve labeled $z$ in the right side of Figure \ref{keylemma1} commutes with the curves labeled $a$ and $b$ in the left side of Figure \ref{keylemma1}. Thus the Dehn twist $z$ commutes with the Dehn twist $ab$ and so we can apply the Key Lemma to $x_0$ and $c$, which are shown in bold in the left of Figure \ref{keylemma1}. Thus the hole encircled by $x_0$ splits and we obtain the relation $zx_1x_0wx_2x_3x_4x_5x_6^2y=abc_1c_2def$, or $x_0x_1x_2x_3x_4x_5x_6^2wyz=abc_1c_2def$, depicted on the right side of Figure \ref{keylemma1}. Notice that in the new relation, we relabeled the boundary parallel curves for convenience and one of the curves that was labeled $x_0$ is now labeled $w$. This relabelling will be done throughout. Once again, it is easy to see that the total space of the Lefschetz fibration described by the monodromy $abc_1c_2def$ is a 2-replacement of the plumbing tree associated to the monodromy $x_0x_1x_2x_3x_4x_5x_6^2wyz$, which is depicted in Figure \ref{plumb1}. Now, inductively assume that the relation $x_0\cdot\cdot\cdot x_{n+3}^2wyz^{n-3}=abc_1\cdot\cdot\cdot c_{n-1}def$ holds, as in the left side of Figure \ref{keylemma2}. Then since $z$ commutes with $ab$, we can apply the Key Lemma to $x_0$ and $c_1$ to obtain the relation $x_0\cdot\cdot\cdot x_{n+3}^2wyz^{n-2}=abc_1\cdot\cdot\cdot c_{n-1}def$. As before, the total space of the Leftschetz fibration described by the monodromy $abc_1\cdot\cdot\cdot c_{n-1}def$ is a 2-replacement of the plumbing tree shown in Figure \ref{plumb2}.

Next, we apply the Key Lemma to $x_{n+2}$ and $e$. To do this, view the $n+3$ punctured disk as an $n+4$ punctured sphere so that the outermost boundary of the disk is just another puncture. In this way, we can view $x_{n+3}$ as a curve around a puncture and $e$ as a curve around the two punctures with boundary parallel curves $x_{n+2}$ and $x_{n+3}$. These are shown in bold on the left side of Figure \ref{keylemma3}. Since $x_{n+4}$ (as labeled on the right side of Figure \ref{keylemma3}) commutes with everything, the Key Lemma applies, yielding the relation $x_0\cdot\cdot\cdot x_{n+1}x_{n+4}x_{n+2}x_{n+3}x_{n+4}^2wyz^{n-2}=abc_1\cdot\cdot\cdot c_{n-1}de_2e_1f$, or $x_0\cdot\cdot\cdot x_{n+3}x_{n+4}^3wyz^{n-2}=abc_1\cdot\cdot\cdot c_{n-1}de_2e_1f$. This relation proves that the linear plumbing depicted in Figure \ref{plumb3} is 2-replaceable. Now, inductively assume that the relation $x_0\cdot\cdot\cdot x_{n+m-2}x_{n+m-1}^{m-2}wyz^{n-2}=abc_1\cdot\cdot\cdot c_{n-1}de_{m-2}\cdot\cdot\cdot e_1f$ holds, as in the left side of Figure \ref{keylemma4}. Notice that each curve $e_i$ encircles all the punctures except the one labeled $i$ for $1\le i\le m-2$, as depicted in Figure \ref{keylemma4}. Also note that, with this labelling, we can write $f=e_0$. We now apply the Key Lemma as we did previously to the bold curves labeled $x_{n+m-2}$ and $e_{m-2}$ to obtain the relation $x_0\cdot\cdot\cdot x_{n+m-1}x_{n+m}^{m-1}wyz^{n-2}=abc_1\cdot\cdot\cdot c_{n-1}de_{m-1}\cdot\cdot\cdot e_1f$. This relation proves that the plumbing tree depicted in Figure \ref{plumb4} is 2-replaceable for all $m\ge3$. Thus we have proved that the plumbing trees of Theorem \ref{thm:trees}a are indeed 2-replaceable.

\begin{figure}[b]
\centering
\begin{subfigure}{\linewidth}
\centering
\includegraphics[scale=.4]{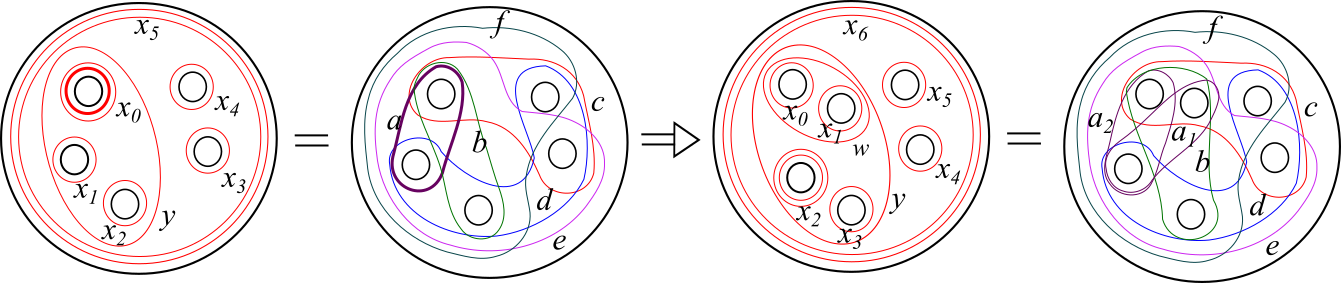}
\caption{$x_0^2x_1x_2x_3x_4x_5^2y=abcdef$ implies $x_0x_1x_2^2x_3x_4x_5x_6^2wy=a_1a_2bcdef$}\label{keylemma5}
\end{subfigure}
\begin{subfigure}{\linewidth}
\centering
\captionsetup{justification=centering}
\includegraphics[scale=.4]{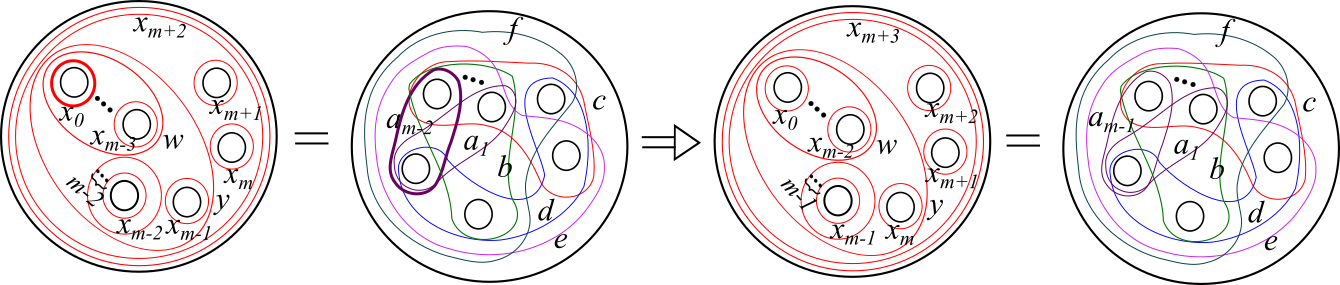}
\caption{$x_0\cdot\cdot\cdot x_{m-3}x_{m-2}^{m-2}x_{m-1}x_mx_{m+1}x_{m+2}^2wy=a_1\cdot\cdot\cdot a_{m-2}bcdef$ implies\\ $x_0\cdot\cdot\cdot x_{m-2}x_{m-1}^{m-1}x_{m}x_{m+1}x_{m+2}x_{m+3}^2wy=a_1\cdot\cdot\cdot a_{m-1}bcdef$}\label{keylemma6} 
\end{subfigure}
\begin{subfigure}{\linewidth}
\centering
\captionsetup{justification=centering}
\includegraphics[scale=.4]{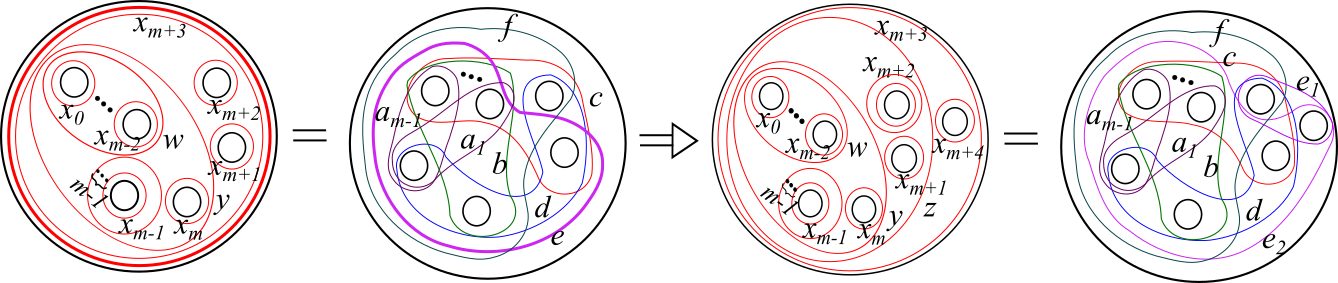}
\caption{$x_0\cdot\cdot\cdot x_{m-2}x_{m-1}^{m-1}x_{m}x_{m+1}x_{m+2}x_{m+3}^2wy=a_1\cdot\cdot\cdot a_{m-1}bcdef$ implies\\ $x_0\cdot\cdot\cdot x_{m-2}x_{m-1}^{m-1}x_{m}x_{m+1}x_{m+2}^2x_{m+3}x_{m+4}wyz=a_1\cdot\cdot\cdot a_{m-1}bcde_1e_2f$}\label{keylemma7}
\end{subfigure}
\begin{subfigure}{\linewidth}
\centering
\captionsetup{justification=centering}
\includegraphics[scale=.4]{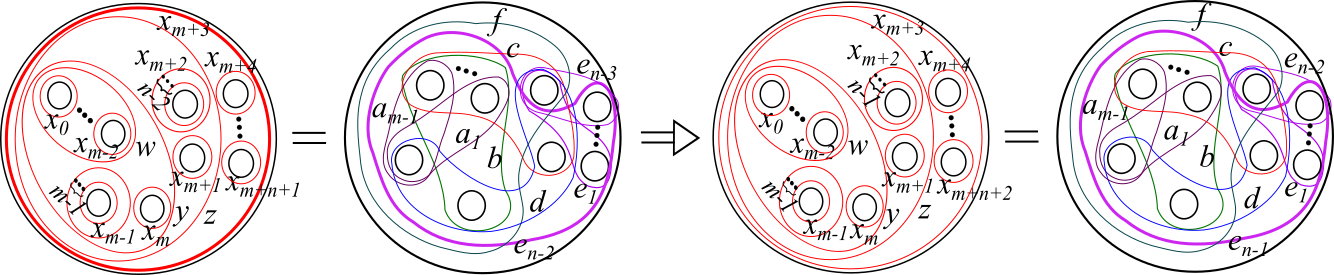}
\caption{$x_0\cdot\cdot\cdot x_{m-2}x_{m-1}^{m-1}x_{m}x_{m+1}x_{m+2}^{n-2}x_{m+3}x_{m+4}wyz=a_1\cdot\cdot\cdot a_{m-1}bcde_1e_2f$ implies\\ $x_0\cdot\cdot\cdot x_{m-2}x_{m-1}^{m-1}x_{m}x_{m+1}x_{m+2}^{n-1}x_{m+3}\cdot\cdot\cdot x_{m+n+2}wyz=a_1\cdot\cdot\cdot a_{m-1}bcde_1\cdot\cdot\cdot e_{n-1}f$}\label{keylemma8} 
\end{subfigure}
\caption{Repeated applications of the Key Lemma to the bold circles}\label{keylemma5-8}
\end{figure}

\begin{figure}[h!]
\centering
\begin{subfigure}[h!]{.4\textwidth}
\centering
\includegraphics[scale=.6]{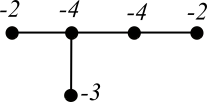}
\caption{Plumbing associated to Figure \ref{keylemma5}}\label{plumb5}
\end{subfigure}
\begin{subfigure}[h!]{.5\textwidth}
\centering
\includegraphics[scale=.6]{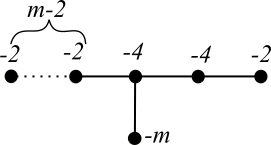}
\caption{Plumbing associated to Figure \ref{keylemma6}}\label{plumb6}
\end{subfigure}
\begin{subfigure}[h!]{.4\textwidth}
\centering
\includegraphics[scale=.6]{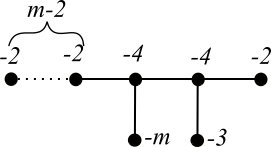}
\caption{Plumbing associated to Figure \ref{keylemma7}}\label{plumb7}
\end{subfigure}
\begin{subfigure}[h!]{.5\textwidth}
\centering
\includegraphics[scale=.6]{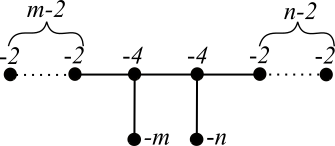}
\caption{Plumbing associated to Figure \ref{keylemma8}}\label{plumb8}
\end{subfigure}
\caption{2-replaceable plumbings associated to the monodromies in Figure \ref{keylemma5-8}}\label{plumb5-8}
\end{figure}

To obtain the family in Theorem \ref{thm:trees}b, we go back to the relation depicted in Figure \ref{relation}, namely $x_0^2x_1x_2x_3x_4x_5^2y=abcdef$. We will apply the Key Lemma to the bold circles $x_0$ and $a$ shown in Figure \ref{keylemma5}.
Since $x_2$, as labeled in the third surface in Figure \ref{keylemma5} commutes with everything, the Key Lemma applies and we obtain the relation $wx_2x_0x_1x_2x_3x_4x_5x_6^2y=a_1a_2bcdef$, or $x_0x_1x_2^2x_3x_4x_5x_6^2wy=a_1a_2bcdef$, as shown in Figure \ref{keylemma5}. Thus, the plumbing tree shown in Figure \ref{plumb5} is 2-replaceable. Inductively assume that $x_0\cdot\cdot\cdot x_{m-3}x_{m-2}^{m-2}x_{m-1}x_mx_{m+1}x_{m+2}^2wy=a_1\cdot\cdot\cdot a_{m-2}bcdef$, as in Figure \ref{keylemma6}. Again, since $x_{m-1}$, as labeled in the third monodromy in Figure \ref{keylemma6}, commutes with everything, we can apply the Key Lemma to obtain the relation $x_0\cdot\cdot\cdot x_{m-2}x_{m-1}^{m-1}x_{m}x_{m+1}x_{m+2}x_{m+3}^2wy=a_1\cdot\cdot\cdot a_{m-1}bcdef$. Thus, the plumbing tree in Figure \ref{plumb6} is 2-replaceable.

Now view the leftmost punctured disk in Figure \ref{relation} as a sphere with six punctures. Then we can arrange the sphere so that the curve labeled $y$ is the equator and the northern and southern hemispheres both have 3 punctures, two of which have one parallel curve and one of which has two parallel curves. In the previous paragraph, we repeatedly applied the Key Lemma to curves in only one of the hemispheres (without involving the equator $y$). Thus we can also apply it to the other hemisphere in the exact same way. We now do this explicitly. In the relation $x_0\cdot\cdot\cdot x_{m-2}x_{m-1}^{m-1}x_{m}x_{m+1}x_{m+2}x_{m+3}^2wy=a_1\cdot\cdot\cdot a_{m-1}bcdef$, consider the bold curves $x_{m+3}$ and $a_{m-1}$ shown in Figure \ref{keylemma6} . We view the latter as a curve containing the two punctures with boundary parallel curves $x_{m+2}$ and $x_{m+3}$. Since $x_{m+2}$ commutes with all other Dehn twists, we can apply the Key Lemma to obtain $x_0\cdot\cdot\cdot x_{m-2}x_{m-1}^{m-1}x_{m}x_{m+1}x_{m+2}x_{m+2}x_{m+3}x_{m+4}zwy=a_1\cdot\cdot\cdot a_{m-1}bcde_1e_2f$, or $x_0\cdot\cdot\cdot x_{m-2}x_{m-1}^{m-1}x_{m}x_{m+1}x_{m+2}^2x_{m+3}x_{m+4}zwy=a_1\cdot\cdot\cdot a_{m-1}bcde_1e_2f$. Thus the plumbing tree in Figure \ref{plumb7} is 2-replaceable. Inductively assume the relation $x_0\cdot\cdot\cdot x_{m-2}x_{m-1}^{m-1}x_{m}x_{m+1}x_{m+2}^{n-2}x_{m+3}\cdot\cdot\cdot x_{m+n+1}zwy=a_1\cdot\cdot\cdot a_{m-1}bcde_1\cdot\cdot\cdot e_{n-2}f$, as in Figure \ref{keylemma8}, holds. Again, since $x_{m+2}$ commutes with everything, we can apply the Key Lemma to obtain the relation $x_0\cdot\cdot\cdot x_{m-2}x_{m-1}^{m-1}x_{m}x_{m+1}x_{m+2}^{n-1}x_{m+3}\cdot\cdot\cdot x_{m+n+2}zwy=a_1\cdot\cdot\cdot a_{m-1}bcde_1\cdot\cdot\cdot e_1\cdot\cdot\cdot e_{n-1}f$. Thus the plumbing tree in Figure \ref{plumb8} is 2-replaceable and so the family of trees in Theorem \ref{thm:trees}b are indeed 2-replaceable.


\section{A symplectic exotic $\mathbb{C}P^2\#6\overline{\mathbb{C}P^2}$}\label{exotic}

In this section we will find the 2-replaceable plumbing tree of Theorem \ref{thm:trees}(a) with $n=9$ and $m=3$ embedded in $\mathbb{C}P^2\#16\overline{\mathbb{C}P^2}$, excise it, and replace it with the 2-replacement constructed in the proof of Theorem \ref{thm:trees}. We will then show that the resulting 4-manifold $X$ is a symplectic exotic $\mathbb{C}P^2\#6\overline{\mathbb{C}P^2}$. We assume the reader is familiar with elliptic fibrations and blowups. See \cite{stipgompf} for details.

\begin{figure}[h!]
\centering
\begin{subfigure}[h!]{.45\textwidth}
\centering
\includegraphics[scale=.6]{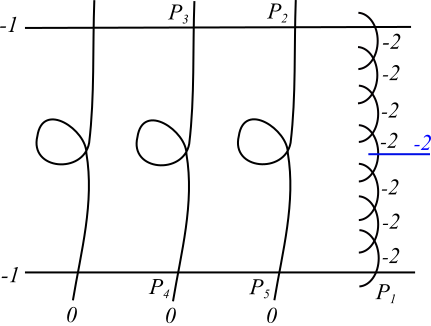}
\caption{A configuration in $E(1)$}\label{e(1)}
\end{subfigure}
\begin{subfigure}[h!]{.45\textwidth}
\centering
\includegraphics[scale=.6]{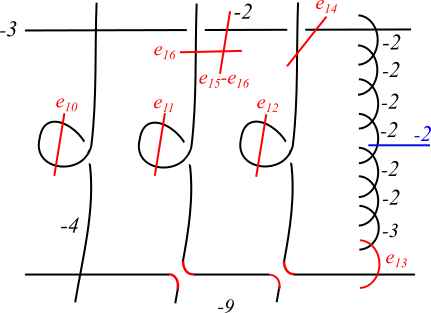}
\caption{Configuration after seven blowups}\label{7blowups}
\end{subfigure}
\caption{Blowing up $E(1)$ seven times}\label{config}
\end{figure}

In \cite{stipsicz}, Stipsicz and Szab\'o showed that there is an elliptic fibration $\mathbb{C}P^2\#9\overline{\mathbb{C}P^2} \rightarrow \mathbb{C}P^1$, called $E(1)$, with three fishtail fibers, two sections, and a singular fiber of type III$^*$ (i.e. an $\tilde{E}_7$ singular fiber) which intersect as in Figure \ref{e(1)}. 

\noindent Starting with this configuration, perform the following moves:
\begin{itemize}
\item blow up the three double points in the fishtail fibers and call the exceptional spheres $e_{10}, e_{11},$ and $e_{12}$;
\item blow up the points $P_1$, $P_2$, and $P_3$ and call the exceptional spheres $e_{13}$, $e_{14}$, and $e_{15}$, respectively;
\item blow up the intersection between $e_{15}$ and image of the adjacent fishtail fiber and call the new exceptional sphere $e_{16}$;
\item and smooth the intersection points $P_4$ and $P_5$.
\end{itemize}
\noindent The resulting configuration is shown in Figure \ref{7blowups}.

\begin{figure}[h!]
\centering
\includegraphics[scale=.6]{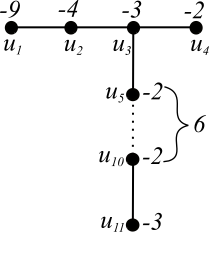}
\caption{The configuration $P$}\label{configC}
\end{figure}

Since we performed seven blowups on $\mathbb{C}P^2\#9\overline{\mathbb{C}P^2}$, this configuration of spheres is embedded in $\mathbb{C}P^2\#16\overline{\mathbb{C}P^2}$. Furthermore, notice that the plumbing $P$ depicted in Figure \ref{configC} is embedded in this configuration. By Theorem \ref{thm:trees}a, this plumbing is 2-replaceable. Let $B$ denote the 2-replacement of $P$ constructed in the proof of Theorem \ref{thm:trees}, let $Z=\mathbb{C}P^2\#16\overline{\mathbb{C}P^2}-\text{int}(P)$ and let $X=Z \cup_{\partial P}B$, where the gluing is by a contactomorphism isotopic to the identity.

\begin{proposition} $X$ is homeomorphic to $\mathbb{C}P^2\#6\overline{\mathbb{C}P^2}$.\label{homeoto}\end{proposition}

\begin{proof} We first prove that $X$ is simply connected. Since $\mathbb{C}P^2\#16\overline{\mathbb{C}P^2}$ is simply connected, the inclusion $\partial P \hookrightarrow Z$ induces a surjection $\pi_1(\partial P)\to\pi_1(Z)$. Furthermore, since $B$ is built out of 0-, 1-, and 2-handles, the inclusion $\partial B=\partial P \hookrightarrow B$ also induces a surjection $\pi_1(\partial P)\to\pi_1(B)$. By the Seifert Van-Kampen theorem, we have $\pi_1(X)=\pi_1(Z)\ast_{\pi_1(\partial P)}\pi_1(B)$. Thus, in the amalgamation, the generators of $\pi_1(Z)$ can be expressed in terms of the generators of $\pi_1(B)$. Therefore, if the generators of $\pi_1(B)$ bound disks in $X$, then $\pi_1(X)$ is trivial. We first prove that $\pi_1(B)$ is cyclic of order $17$ and then show that a particular generator of $\pi_1(B)$ bounds a disk in $X$.

In the proof of Theorem \ref{thm:trees}a, we explicitly described the monodromy of the Lefschetz fibration associated to $B$ (see Figure \ref{keylemma2}). Figure \ref{kirby0} depicts a handlebody diagram of $B$ obtained from this monodromy. For details on how to construct such a diagram, see, for example, \cite{mark}. Each blue unknot has framing $-1$ and, from bottom to top, these unknots correspond to the curves $a,b,c_1,\ldots,c_8,d,e,f$ in the monodromy factorization shown in Figure \ref{keylemma2}. The dotted lines are identified in the trivial way to form an unlink of dotted circles. Let $m_i$ be a meridian around the $i^{th}$ 1-handle of the handlebody diagram of $B$ (shown in Figure \ref{kirby0}), counting left to right. Then $\pi_1(B)$ is generated by $\{m_i\}$ and subject to the following relations (which are given by the 2-handles):
$$m_1\cdot\cdot\cdot  m_9=1,\quad m_1\cdot\cdot\cdot  m_8m_{10}=1, \quad m_1\cdot\cdot\cdot  m_{11}=1,\quad  m_1\cdot\cdot\cdot  m_{10}m_{12}=1,$$ $$m_9\cdot\cdot\cdot m_{12}=1, \text{ and } m_im_{11}m_{12}=1 \text{ for } i=1,\ldots,8. $$

\noindent The last relation shows that $m_1=\cdot\cdot\cdot =m_8$. Call this element $m$. Furthermore, we have $m_{9}=m_{10}=m^{-8}$ and $m_{11}=m_{12}=m^8$. Thus, $1=m_1m_{11}m_{12}=m^{17}$ and so $\pi_1(B)\cong \mathbb{Z}_{17}.$

\begin{figure}[t]
\centering
\begin{subfigure}[h]{.24\linewidth}
\centering
\includegraphics[scale=.5]{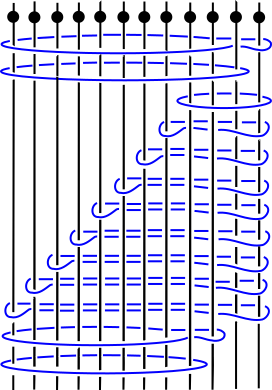}
\caption{A Kirby diagram for $B$}\label{kirby0}
\end{subfigure}
\begin{subfigure}[h]{.24\linewidth}
\centering
\includegraphics[scale=.5]{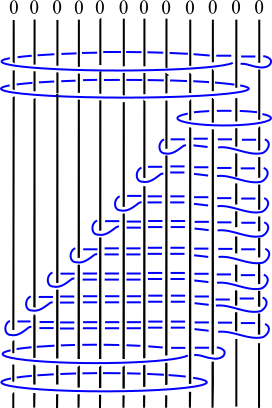}
\caption{A surgery description for $\partial B$}\label{kirby1} 
\end{subfigure}
\begin{subfigure}[h]{.24\linewidth}
\centering
\includegraphics[scale=.5]{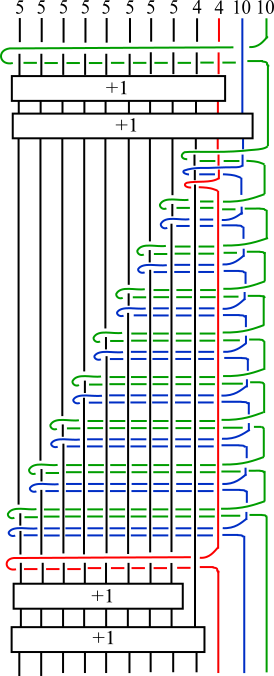}
\caption{Blow down the blue unknots in (B)}\label{kirby2}
\end{subfigure}
\begin{subfigure}[h]{.24\linewidth}
\centering
\includegraphics[scale=.5]{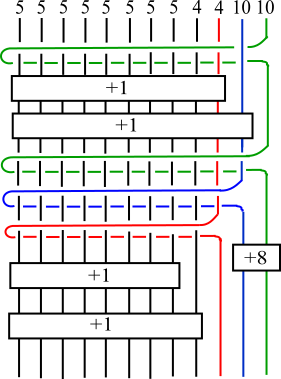}
\caption{Result after isotoping the red, green, and blue strands}\label{kirby3} 
\end{subfigure}
\begin{subfigure}[h]{.24\linewidth}
\centering
\includegraphics[scale=.5]{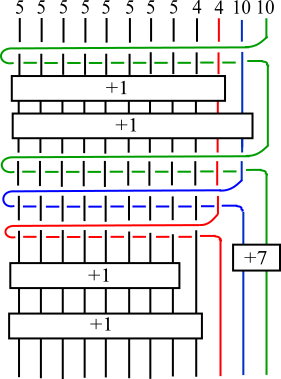}
\caption{Isotope the blue and green strands}\label{kirby4} 
\end{subfigure}
\begin{subfigure}[h]{.24\linewidth}
\centering
\includegraphics[scale=.5]{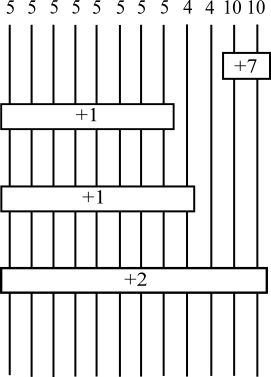}
\caption{Rearrangement of the twists}\label{kirby5}
\end{subfigure}
\begin{subfigure}[h]{.24\linewidth}
\centering
\includegraphics[scale=.5]{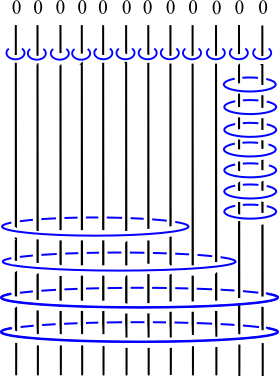}
\caption{Perform blowups}\label{kirby6} 
\end{subfigure}
\begin{subfigure}[h]{.24\linewidth}
\centering
\includegraphics[scale=.5]{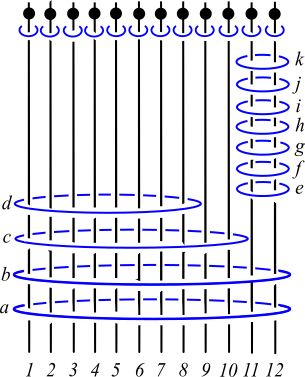}
\caption{A Kirby diagram for $P$}\label{kirby7} 
\end{subfigure}
\caption{Using Kirby calculus to show $B$ and $P$ have the same boundary}\label{kirbycalculus}
\end{figure}

We now use Kirby calculus to move from a handlebody diagram of $B$ to a handlebody diagram of $P$. Start with the handlebody diagram for $B$ depicted in Figure \ref{kirby0} and change the dotted circles to 0-framed unknots to obtain the surgery diagram for $\partial B$ depicted in Figure \ref{kirby1}. Then:
\begin{itemize}
\item blow down all of the blue $-1$ framed unknots to obtain Figure \ref{kirby2}; 
\item isotope the vertical red strand under the strand immediately to its left and pull it leftward;
\item pull the blue and green strands leftward to obtain Figure \ref{kirby3};
\item introduce a positive twist at the top of the blue and green strands and a negative twist at the bottom of the same strands (these twists undo each other) to obtain Figure \ref{kirby4};
\item rearrange the strands to appear as in Figure \ref{kirby5};
\item and perform 23 blowups to obtain Figure \ref{kirby6}. 
\end{itemize}
\noindent Finally, change the 0-framed unknots to dotted circles to obtain the handlebody diagram depicted in Figure \ref{kirby7}. Notice that this is a handlebody diagram for $P$, namely the diagram obtained from the monodromy associated to $P$ in Figure \ref{keylemma2}.

\begin{figure}[h]
\centering
\includegraphics[scale=.55]{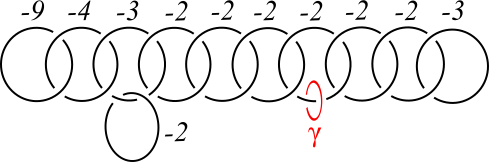}
\caption{Handlebody diagram of $P$ with a meridian $\gamma$}\label{meridian}
\end{figure}

Now consider the obvious handlebody diagram for $P$ depicted in Figure \ref{meridian} (without the red meridian labeled $\gamma$). We can explicitly show that the handlebody diagram in Figure \ref{kirby7} is indeed a handlebody diagram for $P$ via the following moves.

\begin{itemize}
\item Starting with Figure \ref{kirby7}, slide $a$ over $b$, followed by $b$ over $c$, and followed by $c$ over $d$.

\item Slide $d$ over each of the 8 blue unknots at the top of the 1-handles labeled 1-8.

\item Slide $c$ over the blue unknots at the top of the 1-handles labeled 9 and 10.

\item Slide $b$ over $e$, $e$ over $f$, $f$ over $g$, $g$ over $h$, $h$ over $i$, $i$ over $j$, and $j$ over $k$.

\item Slide $k$ over the blue unknots at the top of the 1-handles labeled 11 and 12.

\item Cancel the 1-2 handle pairs to obtain the handlebody diagram in Figure \ref{meridian}.

\end{itemize}

Now consider the red meridian labeled $\gamma$ in Figure \ref{meridian}. By reversing the moves outlined above, in the handlebody diagram in Figure \ref{kirby7}, $\gamma$ links each of the curves labeled $a,b,e,f,g,$ and $h$ exactly once. Changing the dotted circles to 0-framed unknots, we can see $\gamma$ in a surgery diagram of $\partial P$ as in Figure \ref{kirby6withgamma} (imagining that the base point is above the diagram). After blowing down all of the $-1$-framed blue unknots and isotoping $\gamma$, we obtain Figure \ref{kirby5withgamma}. Tracing through the Kirby calculus to obtain the handlebody diagram of $B$ depicted in Figure \ref{kirby0}, it is easy to see that $\gamma$ remains at the bottom of the diagram in the same position as in Figure \ref{kirby5withgamma}. Thus, in $\pi_1(B)$, $\gamma=m_1^2\cdot\cdot\cdot m_{10}^2m_{11}^6m_{12}^6=m^{12}$, which is a generator of $\pi_1(B)\cong\mathbb{Z}_{17}$. 

\begin{figure}
\centering
\begin{subfigure}[h]{.45\linewidth}
\centering
\includegraphics[scale=.55]{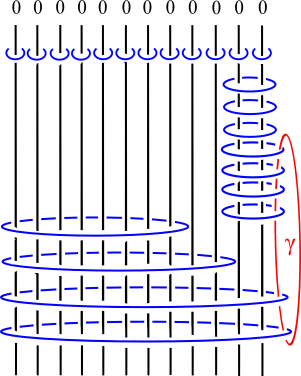}
\caption{$\gamma$ in $\partial P$}\label{kirby6withgamma} 
\end{subfigure}
\begin{subfigure}[h]{.45\linewidth}
\centering
\includegraphics[scale=.55]{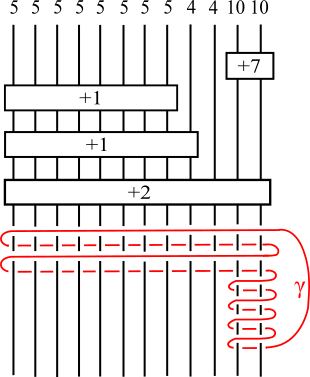}
\caption{$\gamma$ after blowing down}\label{kirby5withgamma} 
\end{subfigure}
\caption{Keeping track of $\gamma$}\label{kirbycalculus}
\end{figure}

Notice, in the original configuration of spheres found in $\mathbb{C}P^2\#16\overline{\mathbb{C}P^2}$ (Figure \ref{7blowups}), $\gamma$ can be identified with the equator of the $-2$-sphere colored in blue that is ``dangling off" the singular fiber of type III*. Thus this meridian bounds a disk (a hemisphere of the blue $-2$-sphere) in $Z$ and thus bounds a disk in $X$. Since $\gamma$ generates $\pi_1(B)$, $X$ is simply connected.

Next, notice that $\chi(X)=\chi(\mathbb{C}P^2\#16\overline{\mathbb{C}P^2})-\chi(P)+\chi(B)=19-12+2=9=\chi(\mathbb{C}P^2\#6\overline{\mathbb{C}P^2})$. Since $\chi(B)=2$ and $b_1(B)=b_3(B)=b_4(B)=0$, we must have $b_2(B)=1$. Since $B$ is negative-definite, the signature of $B$ is $-b_2(B)=-1$ and so $\sigma(X)=\sigma(\mathbb{C}P^2\#16\overline{\mathbb{C}P^2})-\sigma(P)+\sigma(B)=-15-(-11)+(-1)=-5=\sigma(\mathbb{C}P^2\#6\overline{\mathbb{C}P^2})$. Finally, since $-5$ is not divisible by $16$, the intersection forms of $X$ and $\mathbb{C}P^2\#6\overline{\mathbb{C}P^2}$ are both odd. Thus, by Freedman's theorem, $X$ is homeomorphic to $\mathbb{C}P^2\#6\overline{\mathbb{C}P^2}$.\end{proof}

\begin{proposition} $X$ is admits a symplectic structure. \end{proposition}

\begin{proof} Since the spheres of $P$ are complex submanifolds of $\mathbb{C}P^2\#16\overline{\mathbb{C}P^2}$, they are symplectic and intersect positively (c.f. \cite{stipsicz}). By \cite{gompfd3}, these spheres can be made $\omega$-orthogonal by an isotopy through symplectic spheres. By \cite{gaymark}, $X$ admits a symplectic structure.\end{proof}

\begin{proposition} $X$ is not diffeomorphic to $\mathbb{C}P^2\#6\overline{\mathbb{C}P^2}$.\label{exoticprop}\end{proposition}

\begin{proof} Let $h$ denote the canonical generator of $H_2(\mathbb{C}P^2;\mathbb{Z})$ in $H_2(\mathbb{C}P^2\#16\overline{\mathbb{C}P^2};\mathbb{Z})=$\\ $H_2(\mathbb{C}P^2;\mathbb{Z})\bigoplus16H_2(\overline{\mathbb{C}P^2};\mathbb{Z})$ and, with abuse of notation, let $e_i$, for $1\le i\le16$, denote the homology class of the $i^{th}$ exceptional sphere of the $i^{th}$ blowup, which generates the $i^{th}$ copy of $H_2(\overline{\mathbb{C}P^2};\mathbb{Z})$. Consider the configuration of spheres depicted in Figure \ref{e(1)}. Then, as shown in \cite{stipsicz}, the bottom horizontal section has homology class $e_1$, the top horizontal section has homology class $e_9$, each fishtail fiber has homology $3h-\displaystyle\sum_{i=1}^9e_i$, the vertical chain of $-2$-spheres have homology classes (working bottom to top) $h-e_1-e_2-e_3,e_3-e_4,e_4-e_5,e_5-e_6,e_6-e_7,e_7-e_8, e_8-e_9$, and the blue $-2$-sphere has homology class $h-e_3-e_4-e_5$. After performing the seven blowups to obtain the configuration of Figure \ref{7blowups} described earlier, the spheres in our configuration $P$, labeled as in Figure \ref{configC}, have homology classes

$u_1=6h-e_1-\displaystyle\sum_{i=2}^92e_i-2e_{11}-2e_{12}-\sum_{i=13}^{16}e_i, \quad u_2=3h-\displaystyle\sum_{i=1}^9e_i-2e_{10},$
 
$u_3=e_9-e_{14}-e_{15}-e_{16},\quad u_4=e_{15}-e_{16},\quad u_5=e_8-e_9,$

$u_6=e_7-e_8,\quad u_7=e_6-e_7,\quad u_8=e_5-e_6,\quad u_9=e_4-e_5$

$u_{10}=e_3-e_4, \text{ and } u_{11}=h-e_1-e_2-e_3-e_{13}.$\vspace{.3cm}

For quick expositions of Seiberg-Witten invariants when $b^+=1$, see \cite{finsternlectures} and \cite{stipsicz}. It is known that the small perturbation Seiberg-Witten invariant $SW^{\circ}_{\mathbb{C}P^2\#6\overline{\mathbb{C}P^2}}$ is identically 0 because $\mathbb{C}P^2\#6\overline{\mathbb{C}P^2}$ admits a metric of positive scalar curvature. Thus, we must find $\tilde{K}\in H^2(X;\mathbb{Z})$ such that $SW^{\circ}_X(\tilde{K})\neq0$. Let $K$ be the canonical class of $\mathbb{C}P^2\#16\overline{\mathbb{C}P^2}$ associated to the canonical symplectic form $\omega$ on $\mathbb{C}P^2\#16\overline{\mathbb{C}P^2}$. Then $K$ is of the form $K=PD(-3h+\displaystyle\sum_{i=1}^{16}e_i)$ Let $\tilde{K}$ be the canonical class of $X$, induced by the symplectic structure $\tilde{\omega}$ on $X$. Since, by construction, $\omega|_Z=\tilde{\omega}|_Z$, we necessarily have that $K|_Z=\tilde{K}|_Z.$ Furthermore, the dimensions of the Seiberg-Witten moduli spaces associated to $K$ and $\tilde{K}$ are both 0.


By the proof of Corollary 9.4 in \cite{krommrow}, $\partial P$ is an $L$-space. Since $P$ and $B$ are both negative-definite, by Michalogiorgaki's gluing formula in \cite{micha}, $SW^{\circ}_{X,PD(a_2)}(\tilde{K})=SW^{\circ}_{\mathbb{C}P^2\#16\overline{\mathbb{C}P^2},PD(a_1)}(K)$, where $a_1\in H_2(\mathbb{C}P^2\#16\overline{\mathbb{C}P^2};\mathbb{Z})$ and $a_2\in H_2(X;\mathbb{Z})$ such that $a_1|_Z=a_2|_Z$ and $a_1|_P=a_2|_B=0$. Let $$a = \displaystyle 10h-3e_1-2e_2-\sum_{i=3}^93e_i-2e_{10}-e_{11}-2e_{12}-2e_{13}-3e_{14}.$$ Then  $a \cdot u_i=0$ for all $1\le i\le11$ and so $a|_P=0$. Thus $a$ is represented in $Z$ and can also be thought of as a homology class in $H_2(X;\mathbb{Z})$ such that $a|_B=0$. Thus we have $$SW^{\circ}_{X,PD(a)}(\tilde{K})=SW^{\circ}_{\mathbb{C}P^2\#16\overline{\mathbb{C}P^2},PD(a)}(K).$$ 

Since the cohomology class $PD(h)$ gives the chamber that contains the point of positive scalar curvature, $SW^{\circ}_{\mathbb{C}P^2\#16\overline{\mathbb{C}P^2},PD(h)}=0$ (see, e.g. \cite{finsternlectures}). Since $a\cdot a>0$, $h\cdot h>0$, $K\cdot PD(h)=-3<0$, $K\cdot PD(a)=6>0$ and $h\cdot a=10>0$, by the wall crossing formula, we have 
$$SW^{\circ}_{\mathbb{C}P^2\#16\overline{\mathbb{C}P^2},PD(h)}(K)-SW^{\circ}_{\mathbb{C}P^2\#16\overline{\mathbb{C}P^2},\alpha}(K)=(-1)^{1+d(k)/2}$$ and so $SW^{\circ}_{X,\alpha}(\tilde{K})=SW^{\circ}_{\mathbb{C}P^2\#16\overline{\mathbb{C}P^2},\alpha}(K)\neq0$.\end{proof}



\section{Continued fractions}\label{contfracsection}

In this section we outline and prove useful facts about Hirzebruch-Jung continued fractions that will be needed for the proof of Theorem \ref{thm:main}. Given a sequence of integers $(a_1,\ldots,a_n)$ the (Hirzebruch-Jung) continued fraction expansion is given by
\begin{center}
$[a_1,\ldots,a_n]=\displaystyle a_1-\frac{1}{a_2-\displaystyle\frac{1}{\cdot\cdot\cdot-\displaystyle\frac{1}{a_n}}}$
\end{center}\bigskip
If $a_i\ge2$ for all $1\le i\le n$, then this fraction is well-defined and the numerator is greater than the denominator. In fact, for coprime $p>q>0\in\mathbb{Z}$, there exists a unique continued fraction expansion $[a_1,\ldots,a_n]=\frac{p}{q}$, where $a_i\ge2$ for all $1\le i\le n$. To simplify notation, we write $[\ldots,2^{[k]},\ldots]$ instead of $[\ldots,\overbrace{2,\ldots,2}^k,\ldots]$. Moreover, we will often refer to continued fractions as fractions.

We call the continued fraction expansions of $\frac{p}{q}$ and $\frac{p}{p-q}$ \textit{dual} to each other. The following relationship between these two continued fractions is well-known (see, for example, Theorem 7.1 and Lemma 7.2 of \cite{neumann}).

\begin{theorem} Let $n_i\ge 0$ for all $1\le i\le s+1$ and $m_j\ge 0$ for all $1\le j\le s$. If $$\frac{p}{q}=[2^{[n_1]},m_1+3,2^{[n_2]},m_2+3,\ldots,m_s+3,2^{[n_{s+1}]}]$$ then $$\frac{p}{p-q}=[n_1+2,2^{[m_1]},n_2+3,2^{[m_2]},\ldots,n_{s}+3,2^{[m_s]},n_{s+1}+2]$$\label{contfrac}\end{theorem}

\noindent The following corollary follows from Theorem \ref{contfrac}. It will be used throughout the proof of Theorem \ref{thm:main} in Section \ref{mainproof}.

\begin{corollary} If $[m_1,\ldots,m_r]$ has dual $[a_1,\ldots,a_n]$ and $[s_1,\ldots,s_l]$ has dual $[b_1,\ldots,b_k]$, then $[m_1,\ldots,m_r,s_1,\ldots,s_l]$ has dual $[a_1,\ldots,a_{n-1},a_n+b_1-1,b_2,\ldots,b_k]$. Conversely, suppose that $[m_1,\ldots,m_r,s_1,\ldots,s_l]$ has dual $[a_1,\ldots,a_n]$. Then $[m_1,\ldots,m_k]$ and $[s_1,\ldots,s_l]$ have duals of the form $[a_1,\ldots,a_{i-1},a_i']$ and $[a_i'',a_{i+1},\ldots,a_n]$, where $a_i'+a_i''-1=a_i$ and $1\le i\le n$.\label{lem:13}\end{corollary}

\begin{definition} The \textit{buddings} of the fraction $[a_1,\ldots,a_n]$ are the fractions $[a_1+1,a_2,\ldots,a_n,2]$ and $[2,a_1,\ldots,a_{n-1},a_n+1]$. The \textit{debudding} of $[a_1,\ldots,a_n]$ is the reverse operation. (Note: to be able to perform a debudding, we must have either $a_1=2$ and $a_n>2$ or $a_1>2$ and $a_n=2$. For example, the debudding of $[2,a_2,\ldots,a_n]$, where $a_n>2$, is $[a_2,\ldots,a_n-1]$.) Furthermore, by saying $[a_1,\ldots,a_n]$ is a budding of $[a_1',\ldots,a_l']$, we mean that $[a_1,\ldots,a_n]$ can be obtained by a finite sequence of buddings of $[a_1',\ldots,a_l']$ and by saying $[a_1,\ldots,a_n]$ is a debudding of $[a_1',\ldots,a_l']$, we mean that $[a_1,\ldots,a_n]$ can be obtained by a finite sequence of debuddings of $[a_1',\ldots,a_l']$.\end{definition}

Equipped with this definition, the following is a direct consequence of Theorem \ref{contfrac}.

\begin{corollary} If $[a_1,\ldots,a_n]$ has dual $[m_1,\ldots,m_r]$, then the dual of a budding of $[a_1,\ldots,a_n]$ is a budding of $[m_1,\ldots,m_r]$. For example, $[2,a_1,\ldots,a_n+1]$ has dual $[1+m_1,m_2,\ldots,m_r,2]$.\label{cor:2}\end{corollary}

\subsection{Admissible fractions}

In this section, we will consider continued fractions in which all entries are positive and in which each denominator appearing in the fraction is nonzero. Such a fraction is called \textit{admissible}. Note that admissible fractions yield well-defined rational numbers (see, for example, \cite{orlikwag}). In this section, we will consider continued fractions with entries greater than or equal to 1 and so requiring admissibility is important; for example, $[2,1,1]$ is not admissible and is undefined. Moreover, we will consider admissible fractions $[a_1,\ldots,a_n]$ that are equal to 0. The only such continued fraction of length 1 is $[0]$ and if $n\ge 2$, then there must exist an index $i$ such that $a_i=1$ (see, for example, \cite{orlikwag}). 

\begin{definition} Let $[a_1,\ldots,a_n]=0$ be admissible. Then the \textit{blowup before $a_i$} is the fraction $[a_1,\ldots,a_{i-1}+1,1,a_i+1,\ldots,a_n]$ and the \textit{blowup after $a_i$} is the fraction $[a_1,\ldots,a_{i}+1,1,a_{i+1}+1,\ldots,a_n]$. If $a_i=1$, then the \textit{blowdown at $a_i$} is $[a_1,\ldots,a_{i-1}-1,a_{i+1}-1,\ldots,a_n]$. By saying $[a_1,\ldots,a_n]$ is a blowup of $[a_1',\ldots,a_l']$, we mean that $[a_1,\ldots,a_n]$ can be obtained by a finite sequence of blow ups of $[a_1',\ldots,a_l']$. Similarly, $[a_1',\ldots,a_l']$ is a blowdown of $[a_1,\ldots,a_n]$ if it can be obtained by a finite sequence of blow downs of $[a_1,\ldots,a_n]$. \end{definition}

The facts collected in the following proposition are well-known. See, for example, the Appendix of \cite{orlikwag} and Section 2 of \cite{lisca}.

\begin{proposition} If $[a_1,\ldots,a_n]$ is admissible, then:
\begin{enumerate}
\item any blowup or blowdown of $[a_1,\ldots,a_n]$ is also admissible;
\item $[a_n,\ldots,a_1]$ is admissible;
\item $[a_i,a_{i+1},\ldots,a_j]$ is admissible for all $1\le i\le j\le n$;
\item if $[a_1,\ldots,a_n]=0$, then $[a_n,\ldots,a_1]=0$;
\item if $[a_1,\ldots,a_n]=0$, then any blowup or blowdown is also equal to 0; and
\item if $[a_1,\ldots,a_n]=0$, then it can be obtained by a sequence of blowups of $[0].$
\end{enumerate}\label{orlikwagprop}
\end{proposition} 

Note that the only blowup of $[0]$ is $[1,1]$ and the only two blowups of $[1,1]$ are $[1,2,1]$ and $[2,1,2]$. We will consider fractions obtained by sequences of blowups of these two fractions.

\begin{lemma} If $[a_1,\ldots,a_n]=0$ and $a_1=1$ or $a_n=1$, then it is a blowup of $[1,2,1]$.\label{lem:10}\end{lemma}

\begin{proof} We proceed by induction. Let $n=4$. Then the only continued fractions satisfying $[1,a_2,a_3,a_4]=0$ are $[1,3,1,2]$ and $[1,2,2,1]$. By blowing down these continued fractions at the third and fourth entries, respectively, we obtain $[1,2,1]$. Inductively assume that all length $n-1$ (where $n>5$) fractions with $a_1=1$ are blowups of $[1,2,1]$. Let $[a_1,\ldots,a_n]=0$ be a blowup of $[1,2,1]$ and without loss of generality assume $a_1=1$. Then there exists an index $i\neq 1$ such that $a_i=1$. If $i=2$, then since $[1,1,a_3,\ldots,a_n]$ is admissible, so is $[a_n,\ldots,a_3,1,1]$, by Proposition \ref{orlikwagprop}. But $[1,1]=0$, which contradicts admissibility. Thus we may assume $2<i\le n$. By blowing down at $a_i$, we obtain $[1,a_2,\ldots,a_{i-1}-1,a_{i+1}-1,\ldots,a_n]$, which equals 0 and is admissible, by Proposition \ref{orlikwagprop}. By the inductive hypothesis, this fraction is a blowup of $[1,2,1]$. Thus $[a_1,a_2,\ldots,a_n]$ is a blowup of $[1,2,1]$.\end{proof}

\begin{lemma} Let $[a_1,\ldots,a_n]=0$ be a blowup of $[2,1,2]$ that is not a blowup of $[1,2,1]$. Then the buddings of $[a_1,\ldots,a_n]$ are also blowups of $[2,1,2]$ and not of $[1,2,1]$. By Proposition \ref{orlikwagprop}, the buddings are admissible and equal to 0.\label{lem:20}\end{lemma}

\begin{proof} Let $[a_1,\ldots,a_n]$ be as in the statement of the lemma. Then there is a sequence of blowdowns that yields $[2,1,2]$. Performing this sequence of blowdowns to the budding $[2,a_1,\ldots,a_n+1]$, we obtain $[2,2,1,3]$, which is a blowup of $[2,1,2]$. Similarly, $[a_1+1,\ldots,a_n,2]$ is a blowup of $[3,1,2,2]$, which is a blowup of $[2,1,2]$.\end{proof}

\begin{lemma} Let $[a_1,\ldots,a_n]=0$ be a blowup of $[2,1,2]$ that is not a blowup of $[1,2,1]$ and suppose $a_1=2$ and $a_n\ge 3$ (or vica versa). Then the debudding of $[a_1,\ldots,a_n]$ is a blowup of $[2,1,2]$ and not a blowup of $[1,2,1]$. By Proposition \ref{orlikwagprop}, this debudding is admissible and equal to 0.\label{lem:8}\end{lemma}

\begin{proof} We proceed by induction on $n$. First notice that the only blowups of $[2,1,2]$ that are not blowups of $[1,2,1]$ are the fractions $[2,2,1,3]$ and $[3,1,2,2]$. These have one possible debudding each, namely $[2,1,2]$, which is not a blowup of $[1,2,1]$. Inductively assume that the lemma is true for all length $n-1$ continued fractions satisfying the hypotheses. Let $[a_1,\ldots,a_n]=0$ be a blowup of $[2,1,2]$ that is not a blowup of $[1,2,1]$ and suppose without loss of generality that $a_1=2$ and $a_n>2$. Then $a_2\neq 1$, since otherwise, the blowdown at $a_2$ would have first entry equal to 1, which would imply that $[a_1,\ldots,a_n]$ is a blowup of $[1,2,1]$, by Lemma \ref{lem:10}. Let $a_i=1$, where $2<i<n$. By blowing down $[a_1,\ldots,a_n]$ at $a_i$, we obtain the length $n-1$ fraction $[a_1,\ldots,a_{i-1}-1,a_{i+1}-1,\ldots,a_n]=0$, which is not a blowup of $[1,2,1]$. By the inductive hypothesis, the debudding $[a_2,\ldots,a_{i-1}-1,a_{i+1}-1,\ldots,a_n-1]$ is a blowup of $[2,1,2]$ and not a blowup of $[1,2,1]$. Now, by performing a blowup before $a_{i+1}-1$, we obtain $[a_2,\ldots,a_n-1]$, which is the debudding of $[a_1,\ldots,a_n]$. Thus the debudding of $[a_1,\ldots,a_n]$ is a blowup of $[2,1,2]$ and not a blowup of $[1,2,1]$.\end{proof}

By Proposition \ref{orlikwagprop}, Lemma \ref{lem:20}, and Lemma \ref{lem:8}, we will not have to check for the admissibility of any continued fractions throughout the remainder of this article. The final two results are corollaries of the above lemmas. They will be used throughout Section \ref{mainproof}.

\begin{corollary} $[a_1,\ldots,a_n]$ is a budding of $[2,1,2]$ if and only if it is a blowup of $[2,1,2]$ that has exactly one entry that is equal to 1.\label{lem:9}\end{corollary}

\begin{proof} If $[a_1,\ldots,a_n]$ is a budding of $[2,1,2]$, then it necessarily has exactly one entry equal to 1 and by Lemma \ref{lem:20}, it is a blowup of $[2,1,2]$. If $[a_1,\ldots,a_n]$ is a blowup of $[2,1,2]$ with exactly one entry equal to 1, then either $a_1=2$ and $a_n\ge 3$ or vica versa. Thus by repeated applications of Lemma \ref{lem:8}, it is clear that $[a_1,\ldots,a_n]$ is a budding of $[2,1,2]$.\end{proof}

\begin{corollary} $[b_1,\ldots,b_k]$ is a budding of $[4]$ if and only if $[b_1,\ldots,b_k]$ has dual of the form $[a_1,\ldots,a_i+1,\ldots,a_n]$, where $a_i=1$, $1<i<n$, and $[a_1,\ldots,a_i,\ldots,a_n]$ is a blowup of $[2,1,2]$ with exactly one entry that is 1.\label{4}\end{corollary}

\begin{proof} Since the dual of $[4]$ is $[2,2,2]$, this follows from Corollaries \ref{cor:2} and \ref{lem:9}.\end{proof}


\section{Proof of Theorem \ref{thm:main}}\label{mainproof}
\subsection{Lisca's classification of symplectic fillings of $(L(p,q),\xi_{st})$}\label{liscaclass}

Let $p>q>0\in\mathbb{Z}$ be coprime. In \cite{lisca}, Lisca classified all minimal weak symplectic fillings of $(L(p,q),\xi_{st})$, where $\xi_{st}$ is the standard tight contact structure on $L(p,q)$ inherited from the unique tight contact structure on $S^3$. It is known that every weak symplectic filling of a rational homology sphere can be modified into a strong symplectic filling (\cite{ohtaono}). Thus Lisca's classification is of strong symplectic fillings of $(L(p,q),\xi_{st})$. 

\begin{figure}[h!]
\centering
\includegraphics[scale=.8]{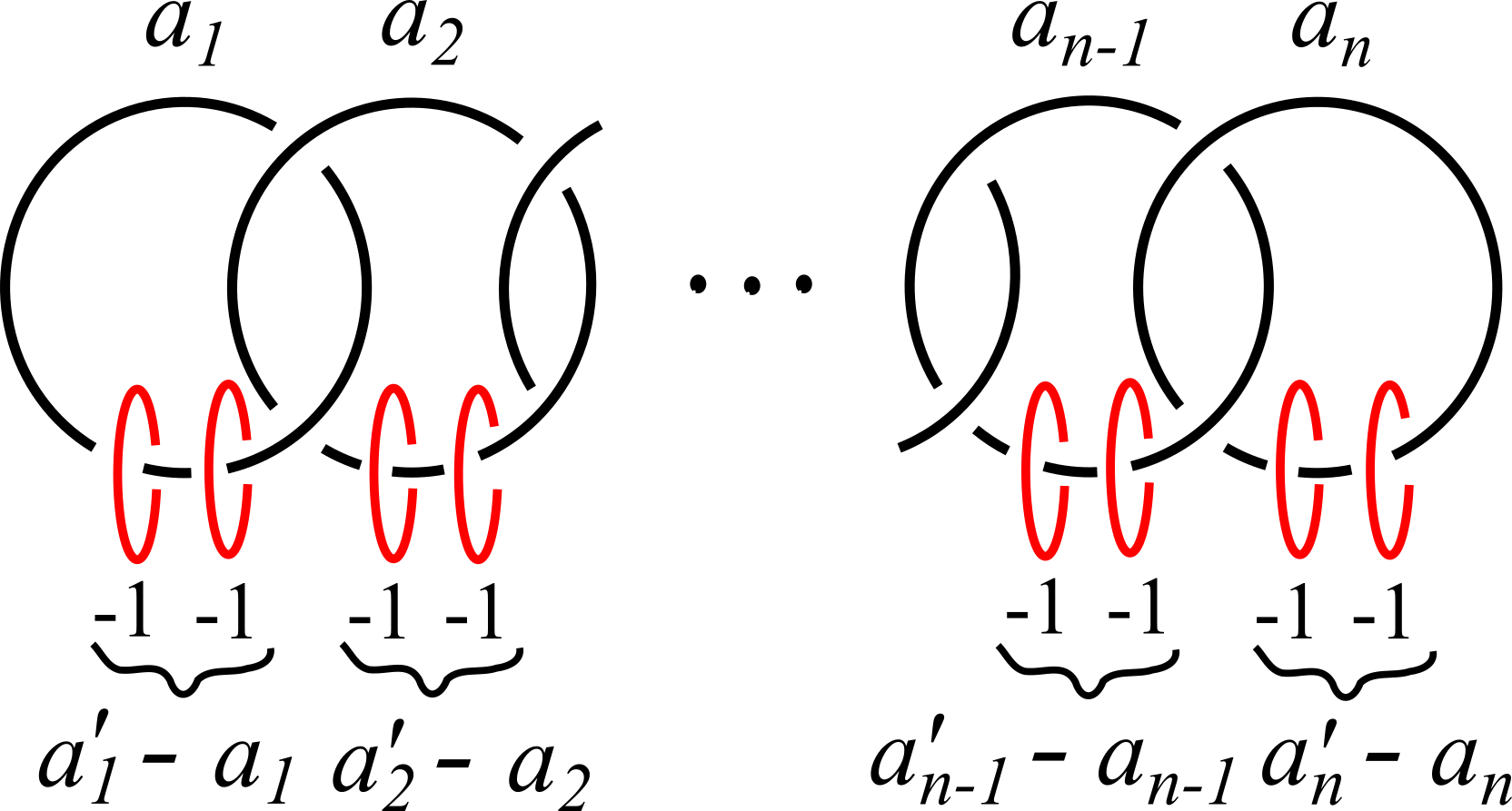}
\caption{A cobordism $W$ from $S^1\times S^2$ to $L(p,q)$, where $\frac{p}{p-q}=[a_1',\ldots,a_n']$ and $[a_1,\ldots,a_n]=0$ is admissible. Any weak (or strong) symplectic filling of $(L(p,q),\xi_{st})$ is obtained by gluing such a cobordism to $S^1\times D^3$.}\label{figure1}
\end{figure}

In particular, Lisca proved that any weak (or strong) symplectic filling of $(L(p,q),\xi_{st})$, where $\frac{p}{p-q}=[a_1',\ldots,a_n']$, is orientation preserving diffeomorphic to a 4-manifold of the form $(S^1\times D^3)\cup W$, where $W$ is the 2-handle cobordism from $S^1\times S^2$ to $L(p,q)$ depicted in Figure \ref{figure1}. In the figure, $[a_1,\ldots,a_n]=0$ is an admissible fraction with $a_i\le a_i'$ for all $1\le i\le n$ and the horizontal chain of framed unknots is a surgery diagram of $S^1\times S^2$. After thickening $S^1\times S^2$ to $S^1\times S^2\times I$, we can attach 2-handles along the $-1$-framed red unknots shown in the diagram to obtain $W$. The angular brackets around the surgery coefficients are simply meant to distinguish the surgery description of $S^1\times S^2$ from the 4-dimensional 2-handles. By gluing $W$ to $S^1\times D^3$, we obtain a filling of $L(p,q)$ with Euler characteristic $\sum_{i=1}^n (a_i'-a_i).$ 

\begin{figure}[h!]
\centering
\includegraphics[scale=.6]{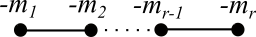}
\caption{Linear plumbing bounded by $L(p,q)$, where $\frac{p}{q}=[m_1,\ldots,m_r]$.}\label{linearfigure}
\end{figure}

From this description, it is easy to see that the first and third Betti numbers of these fillings are 0. Moreover, since $(L(p,q),\xi_{st})$ is supported by a planar open book, every filling is negative-definite (\cite{etnyre}). Now suppose $\frac{p}{p-q}=[a_1',\ldots,a_n']$ has dual $\frac{p}{q}=[m_1,\ldots,m_r]$. Then $L(p,q)$ can be obtained by performing Dehn surgery along a chain of unlinks with surgery coefficients $(-m_1,\ldots, -m_r)$. It follows that $L(p,q)$ bounds the negative-definite linear plumbing $P$ with weights $(-m_1,\ldots,-m_r)$ depicted in Figure \ref{linearfigure}, which is known to be a strong symplectic filling of $(L(p,q),\xi_{st})$. Thus $P$ is $\sum_{i=1}^n (a_i'-a_i)$-replaceable.

\begin{remark} Since $[a_1,\ldots,a_n]=0$ is admissible, we can perform a sequence of blow downs to obtain the continued fraction $[0]$. To see an explicit handlebody diagram of the filling, we can perform the same sequence of blow downs to the horizontal chain of unknots in Figure \ref{figure1} to obtain an unknot with framing $\langle 0\rangle$. In this new diagram, the $\langle0\rangle$-framed unknot represents $S^1\times S^2$, which still bounds $S^1\times D^3$ in the filling. Thus by changing the $\langle0\rangle$-framed unknot to a dotted circle, we obtain an honest handlebody diagram of the filling. Under these moves, the $-1$-framed red unknots become a complicated link with very negative framings.\end{remark}

Lisca's classification shows that there is a one-to-one correspondence between continued fractions of the form $\frac{p}{p-q}=[a_1+k_1,\ldots,a_n+k_n]$, where $[a_1,\ldots,a_n]=0$ is admissible, and (Euler characteristic $\sum_{i=1}^n k_i$) symplectic fillings of $(L(p,q),\xi_{st})$. For example, continued fractions of the form $[a_1,\ldots,a_i+1,\ldots,a_n]=\frac{p}{p-q}$, where $[a_1,\ldots,a_i,\ldots,a_n]=0$ is admissible and $a_i=1$ is the only entry equal to 1, correspond to Euler characteristic 1 symplectic fillings of $(L(p,q),\xi_{st})$. By Corollary \ref{4}, the classification of 1-replaceable linear plumbings (in terms of the weights of the plumbing graph) is immediate.

\begin{corollary} A linear plumbing with weights $(-m_1,\ldots,-m_r)$ is 1-replaceable if and only if $[m_1,\ldots,m_r]$ is a budding of $[4]$.\label{1repcor}\end{corollary}

To prove Theorem \ref{thm:main}, we will classify the continued fractions $\frac{p}{q}$ whose dual continued fractions correspond to Euler characteristic 2 symplectic fillings. To this end, we now assume that every continued fraction $[a_1,\ldots,a_n]=0$ is admissible with at most two entries that are equal to 1. Suppose $a_i=a_j=1$. If $i\neq j$, then $[a_1,\ldots,a_i+1,\ldots,a_j+1,\ldots,a_n]=\frac{p}{p-q}$ corresponds to an Euler characteristic 2 symplectic filling of $(L(p,q),\xi_{st})$ and if $i=j$, then $[a_1,\ldots,a_i+2,\ldots,a_n]=\frac{p_i}{p_i-q_i}$ and $[a_1,\ldots, a_i+1,\ldots,a_t+1,\ldots,a_n]=\frac{p_t}{p_t-q_t},$ where $t\neq i$, correspond to Euler characteristic 2 symplectic fillings of $(L(p_t,q_t),\xi_{st})$ for all $1\le t\le n$. The following corollary summarizes this discussion.

\begin{corollary} The linear plumbing with weights $(-m_1,\ldots,-m_r)$ is 2-replaceable if and only if the dual of $[m_1,\ldots,m_r]$ is of the form:
\begin{itemize}[leftmargin=.9cm]
\item $[a_1,\ldots,a_i+1,\ldots,a_j+1,\ldots,a_n]$, where $[a_1,\ldots,a_i,\ldots,a_j,\ldots,a_n]=0$ has exactly two entries equal to 1, namely $a_i$ and $a_j$;
\item $[a_1,\ldots,a_i+1,\ldots,a_t+1,\ldots,a_n]$, where $t\neq i$ and $[a_1,\ldots,a_i,\ldots,a_j,\ldots,a_n]=0$ has exactly one entry equal to 1, namely $a_i$; or
\item $[a_1,\ldots,a_i+2,\ldots,a_n]$, where $[a_1,\ldots,a_i,\ldots,a_n]=0$ has exactly one entry equal to 1, namely $a_i$.
\end{itemize}
\label{summary}
\end{corollary}

\noindent To prove Theorem \ref{thm:main}, we will explore the continued fractions listed in Corollary \ref{summary}.

\subsection{Proof of Theorem \ref{thm:main}}
For convenience, we recall Theorem \ref{thm:main}. \\

\noindent \textbf{Theorem \ref{thm:main}:}\textit{ Let $(-b_1,\ldots,-b_k)$ and $(-c_1,\ldots,-c_l)$ be obtained by sequences of buddings of $(-4)$ and let $z\ge2$ be any integer. Then a linear plumbing is 2-replaceable if and only if it is either of the form:}
\vspace{.1cm}
\begin{enumerate}[label=(\alph*)]
\item \includegraphics[scale=.6]{linearplumbing1} \qquad\qquad\qquad\qquad  for $k,l\ge0$
\end{enumerate}
\textit{or is obtained by a sequence of buddings of one of the linear plumbings of the form}:
\begin{enumerate}[resume, label=(\alph*)]
\item \includegraphics[scale=.6]{linearplumbing2}\quad (or \quad \includegraphics[scale=.6]{linearplumbing3}) \qquad\quad for $k\ge0$.
\vspace{.1cm}
\item \includegraphics[scale=.6]{linearplumbing4}
\vspace{.1cm}
\item \includegraphics[scale=.6]{linearplumbing5} \qquad\qquad\quad for $k,l\ge1$
\end{enumerate}
\vspace{.3cm}

First, we first show that all plumbings listed in the theorem are indeed 2-replaceable. This is clear for the linear plumbings of type (a), since the ``subplumbings" on either side of the $-z$-disk bundle can be symplectically rationally blown down, revealing an Euler characteristic 2 symplectic 4-manifold. Buddings of plumbings of type (b), (c), and (d) can also be seen to be 2-replaceable via Kirby calculus, blowups, and rational blowdowns. Instead of working through these details, however, we will apply Corollary \ref{summary}.\\


\textbf{Type (b)}: Suppose $[b_1,\ldots,b_k]$ is a budding of $[4]$. By Corollary \ref{4},  $[b_1,\ldots,b_k]$ has dual of the form $[a_1,\ldots,a_i+1,\ldots,a_n]$, where $[a_1,\ldots,a_i,\ldots,a_n]=0$ has exactly one entry that is equal to 1, namely $a_i$, where $i\neq 1,n$. By Theorem \ref{contfrac}, $[b_1,\ldots,b_k,2]$ has dual $[a_1,\ldots,a_i+1,\ldots,a_n+1]$ and so by Corollary \ref{summary} the plumbing with weights $(-b_1,\ldots,-b_k,-2)$ is 2-replaceable. Now let $[m_1,\ldots,m_r]$ be a budding of $[b_1,\ldots,b_k,2]$. Then, by Corollary \ref{cor:2}, it has dual that is a budding of $[a_1,\ldots,a_i+1,\ldots,a_n+1]$ and is of the form $[\ldots,a_1+t,a_2,\ldots,a_i+1,\ldots,a_n+1+s,\ldots]$, for some $s,t\ge0$. By Corollary \ref{lem:9}, $[\ldots,a_1+t,a_2,\ldots,a_i,\ldots,a_n+s,\ldots]=0$ and it is a blowup of $[2,1,2]$ with exactly one entry equal to 1, namely $a_i$. Thus by Corollary \ref{summary}, the plumbing with weights $(-m_1,\ldots,-m_r)$ is 2-replaceable. Similarly, any budding of the linear plumbing with weights $(-2,-b_1,\dots,-b_k)$ is 2-replaceable.\qed\\

\textbf{Type (c)}: $[3,3]$ has dual $[2,3,2]$. Since $[2,1,2]=0$, by Corollary \ref{summary}, the plumbing with weights $(-3,-3)$ is 2-replaceable. Let $[m_1,\ldots,m_r]$ be a budding of $[3,3]$. By applying Corollary \ref{cor:2}, Corollary \ref{lem:9}, and Corollary \ref{summary} as in the proof of ``Type (b)," it is clear that the plumbing with weights $(-m_1,\ldots,-m_r)$ is 2-replaceable.\qed\\

\textbf{Type (d)}: Let $[b_1,\ldots,b_k]$ and $[c_1,\ldots,c_l]$ be buddings of $[4]$. Then they have respective duals of the form $[a_1,\ldots,a_i+1,\ldots,a_n]$ and $[a_1',\ldots,a_j'+1,\ldots,a_m']$, where $[a_1,\ldots,a_i,\ldots,a_n]=[a_1',\ldots,a_j',\ldots,a_m']=0$ and $a_i=a_j'=1$ are the only entries equal to 1. By Corollary \ref{lem:13}, $[b_1,\ldots,b_k,c_1,\ldots,c_l]$ has dual $[a_1,\ldots,a_i+1,\ldots,a_{n-1},a_n+a_1'-1,a_2',\ldots,a_j'+1,\ldots,a_m']$ and so, by Theorem \ref{contfrac}, $[2,b_1,\ldots,b_k,c_1,\ldots,c_l,2]$ has dual $[a_1+1,\ldots,a_i+1,\ldots,a_{n-1},a_n+a_1'-1,a_2',\ldots,a_j'+1,\ldots,a_m'+1]$. We now claim that $[a_1+1,\ldots,a_i,\ldots,a_{n-1},a_n+a_1'-1,a_2',\ldots,a_j',\ldots,a_m'+1]=0$. By Lemma \ref{lem:20}, $[2,a_1',\ldots,a_j',\ldots,a_m'+1]=0$. It follows that $[a_1',\ldots,a_j',\ldots,a_m'+1]=\frac{1}{2}$ and so $[a_n+a_1'-1,a_2',\ldots,a_j',\ldots,a_m'+1]=a_n-1+[a_1',\ldots,a_j',\ldots,a_m']=a_n-\frac{1}{2}$. Once again, by Lemma \ref{lem:20}, $[a_1+1,\ldots,a_i,\ldots,a_{n-1},a_n-\frac{1}{2}]=[a_1+1,\ldots,a_i,\ldots,a_{n-1},a_n,2]=0$, which implies that $[a_1+1,\ldots,a_i,\ldots,a_{n-1},a_n+a_1'-1,a_2',\ldots,a_j',\ldots,a_m'+1]=[a_1+1,\ldots,a_i,\ldots,a_{n-1},a_n-\frac{1}{2}]=0$. Therefore, by Corollary \ref{summary}, the plumbing with weights $(-2,-b_1,\ldots,-b_k,-c_1,\ldots,-c_l,-2)$ is 2-replaceable. Finally suppose $[m_1,\ldots,m_r]$ is a budding of $[2,b_1,\ldots,b_k,c_1,\ldots,c_l,2]$. By applying Corollary \ref{cor:2}, Corollary \ref{lem:9}, and Corollary \ref{summary} as in the proof of ``Type (b)," it is clear that the plumbing with weights $(-m_1,\ldots,-m_r)$ is 2-replaceable.\qed\\

We have shown that all of the linear plumbings listed in Theorem \ref{thm:main} are indeed 2-replaceable. Next we show that these are the only 2-replaceable linear plumbings. To do this, we consider the continued fractions listed in Corollary \ref{summary} and show that their dual continued fractions \textit{correspond} to the linear plumbings listed in Theorem \ref{thm:main}. That is, we will show that if such a dual is given by $[m_1,\ldots,m_r]$, then the plumbing with weights $(-m_1,\ldots,-m_r)$ is one of the plumbings listed in Theorem \ref{thm:main}.

The first two cases are the continued fractions $[0]$ and $[1,1]$. Adding 2 to the only entry of $[0]$ yields the fraction $[2]=\frac{2}{1}$, which has dual fraction $[2]=\frac{2}{1}$. This corresponds to the $-2$-disk bundle over $S^2$. Similarly, adding 1 to each entry in $[1,1]$ yields the fraction $[2,2]=\frac{3}{2}$, which has dual fraction $[3]=\frac{3}{1}$. This corresponds to the $-3$-disk bundle over $S^2$. These two plumbings are \textit{trivially} 2-replaceable and fall under Theorem \ref{thm:main}(a). Now, there are two blowups of $[1,1]$, namely $[1,2,1]$ and $[2,1,2]$. Since any admissible fraction $[a_1,\ldots,a_n]=0$ with $n\ge 3$ is necessarily a blowup of one of these two fractions (by Proposition \ref{orlikwagprop}), we will consider the following types of continued fractions.

\begin{enumerate}
	\item[1.] Blowups of $[1,2,1]$ with two entries equal to 1 and with either:
	\begin{enumerate}[label=(\roman*)]
		\item $a_1=a_n=1$;
		\item $a_1=1$ and $a_n\neq1$ (or vica versa); or
		\item $a_1,a_n\neq1$;
	\end{enumerate}	
	\item[2.] Blowups of $[2,1,2]$ that are not blowups of $[1,2,1]$ and with either:
	\begin{enumerate}[label=(\roman*)]
		\item exactly one entry that is equal to 1; or
		\item exactly two entries that are equal to 1.
	\end{enumerate}
\end{enumerate}	
\vspace{.1cm}

\noindent For these five types of continued fractions, we can proceed as in the cases of the continued fractions $[0]$ and $[1,1]$. That is, for continued fractions of type 2(i), we can add 2 to the only entry that is equal to 1 (e.g. $[3,1,2,2]$ becomes $[3,3,2,2]$), and for continued fractions of type 1(i),1(ii),1(iii), or 2(ii), we can add 1 to each entry that is equal to 1 (e.g. $[1,3,1,2]$ becomes $[2,3,2,2]$). We claim that the resulting continued fractions have dual continued fractions corresponding to the plumbings listed in Theorem \ref{thm:main}. We will prove this claim case by case in the following five lemmas. Once completed, we will have finished the proof of Theorem \ref{thm:main}.

\begin{lemma}[Type 1(i)] Suppose $[a_1,\ldots,a_n]=0$ is a blowup of $[1,2,1]$ with $a_1=a_n=1$. Then $[a_1+1,a_2,\ldots,a_{n-1},a_n+1]$ has dual of the form $[n+1]$, which corresponds to a plumbing in Theorem \ref{thm:main}(a). \label{lemforthm0}\end{lemma}

\begin{proof} We first claim that $[a_1,\ldots,a_n]$ is of the form $[1,2^{[n-2]},1]$. Otherwise, if there exists an index $t$ such that $a_t\ge 3$, then we can repeatedly blow down the first and last entries until we obtain a continued fraction that is not equal to 0, which contradicts Proposition \ref{orlikwagprop}. Thus $[a_1,\ldots,a_n]$ is of the form $[1,2^{[n-2]},1]$. Adding 1 to the first and last entries yields $[2^{[n]}]=\frac{n+1}{n}$, which has dual $[n+1]=\frac{n+1}{1}$.\end{proof}

\begin{lemma}[Type 1(ii)] Suppose $[a_1,\ldots,a_n]=0$ is a blowup of $[1,2,1]$ with $a_i=1$, where $1<i<n$ and either $a_1=1$ or $a_n=1$. If $a_1=1$, then $[a_1+1,\ldots,a_i+1,\ldots,a_n]$ has dual of the form $[z,c_1,\ldots,c_l]$, where $[c_1,\ldots,c_l]$ is a budding of $[4]$ and $z\ge3$. If $a_n=1$ then $[a_1,\ldots,a_i+1,\ldots,a_n+1]$ has dual of the form $[b_1,\ldots,b_k,z]$, where $[b_1,\ldots,b_k]$ is a budding of $[4]$ and $z\ge3$. Moreover, these duals correspond to plumbings in Theorem \ref{thm:main}(a). \label{lemforthm1}\end{lemma}

\begin{proof} Let $[a_1,\ldots,a_i,\ldots,a_n]=0$ be a blowup of $[1,2,1]$ with $a_1=a_i=1$, where $i\neq n$. First notice that there exists $2\le t<i$ such that $a_t\ge3$. Otherwise, $[a_1,\ldots,a_n]=[1,2,\ldots,2,1,a_{i+1},\ldots,a_n]=0$ and so $[a_n,\ldots,a_{i+1},1,2,\ldots,2,1]=0$, by Proposition \ref{orlikwagprop}. But $[1,2,\ldots,2,1]=0$, which implies that $[a_n,\ldots,a_{i+1},1,2,\ldots,2,1]$ is undefined. Thus $a_2=a_3=\ldots=a_{t-1}=2$ and $a_t\ge3$, where $t<i$. Blow down  the fraction repeatedly at the first entry until we obtain $[a_t-1,\ldots,a_i,\ldots,a_n]=0$. This fraction has exactly one entry that is equal to 1, namely $a_i$. Thus it is a blowup of $[2,1,2]$ and so by Corollary \ref{4}, the dual of $[a_t-1,\ldots,a_i+1,\ldots,a_n]$ is a budding of $[4]$. Denote this dual by $[c_1,\ldots,c_l]$. By Theorem \ref{contfrac}, $[a_t,\ldots,a_i+1,\ldots,a_n]$ has dual $[2,c_1,\ldots,c_l]$. By applying Theorem \ref{contfrac} again, we have that $[a_1+1,\ldots,a_t,\ldots,a_i+1,\ldots,a_n]=[2^{[t-1]},a_t,\ldots,a_i+1,\ldots,a_n]$ has dual $[t+1,c_1,\ldots,c_l]$, where $t+1\ge3$. Setting $z=t+1$, we have the result.

Now suppose $\displaystyle[a_1,\ldots,a_n]=0$ is a blowup of $[1,2,1]$ with $a_i=a_n=1$, where $i\neq 1$. Then by Proposition \ref{orlikwagprop}, $[a_n,\ldots,a_1]=0$ and so, by the first part of the lemma, $[a_n+1,\ldots,a_i+1,\ldots,a_1]$ has dual of the form $[z,b_k,\ldots,b_1]$, where $[b_k,\ldots,b_1]$ is a budding of $[4]$. Thus, by Theorem \ref{contfrac}, the dual of $[a_1,\ldots,a_i+1,\ldots,a_n+1]$ is $[b_1,\ldots,b_k,z]$. Since $[b_1,\ldots,b_k]$ is a budding of $[4]$ and $z\ge 3$, we have the result.\end{proof}

\begin{lemma}[Type 1(iii)] Let $[a_1,\ldots,a_n]=0$ be a blowup of $[1,2,1]$ such that $a_i=a_j=1$, where $1<i<j<n$. Then $[a_1,\ldots,a_i+1,\ldots,a_j+1,\ldots,a_n]$ has dual $[b_1,\ldots,b_k,z,c_1,\ldots,c_l]$, where $[b_1,\ldots,b_k]$ and $[c_1,\ldots,c_l]$ are buddings of $[4]$ and $z\ge2$. Moreover, this dual corresponds to a plumbing in Theorem \ref{thm:main}(a).\label{lemforthm2}\end{lemma}

\begin{proof} Let $[a_1,\ldots,a_n]=0$ be as in the statement of the lemma. Since $[a_1,\ldots,a_n]$ is a blowup of $[1,2,1]$, we can repeatedly blow down at the first occurrence of 1 until the first entry is equal to 1. This blowdown is of the form $[1,a_r-m,\ldots,a_j,\ldots,a_n]=0$, where $r<j$ and $m\le a_r-2$. Moreover, assume that we performed the minimal number of blowdowns to obtain such a fraction. By Lemma \ref{lemforthm1}, $[2,a_r-m,\ldots,a_j+1,\ldots,a_n]$ has dual of the form $[y,c_1,\ldots,c_l]$, where $y\ge3$ and $[c_1,\ldots,c_l]$ is a budding of $[4]$. To recover the original fraction, we now perform blowups. The first blowup must occur after the first entry, since otherwise, we would obtain $[1,2,a_r-m,\ldots,a_n]$, which contradicts the minimality assumption. Thus the first blowup yields $[2,1,a_r-m+1,\ldots,a_n]=0$. By Theorem \ref{contfrac}, $[a_r-m,\ldots,a_j+1,\ldots,a_n]$ has dual $[y-1,c_1,\ldots,c_l]$. Thus $[a_r-m+1,\ldots,a_j+1,\ldots,a_n]$ has dual $[2,y-1,c_1,\ldots,c_l]$ and so $[2,2,a_r-m+1,\ldots,a_j+1,\ldots,a_n]$ has dual $[4,y-1,c_1,\ldots,c_l]$. Writing $[2,2,a_r-m+1,\ldots,a_j+1,\ldots,a_n]$ as $[\underline{2,2,2}+(a_r-m-1),\ldots,a_j+1,\ldots,a_n]$, we can view the underlined portion $[2,2,2]$ as a ``subfraction" of $[2,2,a_r-m+1,\ldots,a_j+1,\ldots,a_n]$. Recall that $[2,2,2]$ has dual $[4]$.

Now blow up $[\underline{2,1,2}+a_r-m-1,\ldots,a_n]=0$ repeatedly before or after the first occurrence of 1 to recover the original fraction $[a_1,\ldots,a_i,\ldots,a_r,\ldots,a_j,\ldots,a_n]=[a_1,\ldots,a_i,\ldots,\\ m+1+(a_r-m-1),\ldots,a_j,\ldots,a_n]=0$. By doing this, we also end up blowing up the subfraction $[2,1,2]$ before or after the only entry that is 1 to obtain $[a_1,\ldots,a_i,\ldots,a_{r-1},m+1]$, which is a subfraction of $[a_1,\ldots,a_i,\ldots,a_r,\ldots,a_j,\ldots,a_n]$. Thus by Corollary \ref{4}, $[a_1,\ldots,a_i+1,\ldots,a_{r-1},m+1]$ has dual that is a budding of $[4]$; denote this budding by $[b_1,\ldots,b_k]$. Since $[a_r-m,\ldots,a_j+1,\ldots,a_n]$ has dual $[y-1,c_1,\ldots,c_l]$, by Lemma \ref{lem:13}, $[a_1,\ldots,a_i+1,\ldots,a_r,\ldots,a_j+1,\ldots,a_n]$ has dual $[b_1,\ldots,b_k,y-1,c_1,\ldots,c_l]$, where $y-1\ge2$. Setting $z=y-1$, we have obtained the result.\end{proof}

\begin{lemma}[Type 2(i)] Let $[a_1,\ldots,a_n]=0$ be a blowup of $[2,1,2]$, where $a_i=1$ is the only entry that is equal to 1 (where $1<i<n$). If $t\neq i$, then the dual fraction of $[a_1,\ldots,a_i+1,\ldots,a_t+1,\ldots,a_n]$ is a budding of either $[2,b_1,\ldots,b_k]$ or $[b_1,\ldots,b_k,2]$, where $[b_1,\ldots,b_k]$ is a budding of $[4]$. Moreover, these buddings correspond to plumbings in Theorem \ref{thm:main}(b). If $t=i$, then the dual fraction of $[a_1,\ldots,a_i+2,\ldots,a_n]$ is a budding of $[3,3]$, which corresponds to a plumbing in Theorem \ref{thm:main}(c).\label{lemforthm3}\end{lemma}

\begin{proof} Suppose $t<i$. Since $[a_1,\ldots,a_n]=0$ is a blowup of $[2,1,2]$ and $a_i$ is the only entry that is 1, it is a budding of $[2,1,2]$, by Corollary \ref{lem:9}. Furthermore, we can view $a_i$ as the image of 1 after performing the sequence of buddings of $[2,1,2]$ to obtain $[a_1,\ldots,a_n]$. Thus we can perform debuddings until we obtain $[a_t,\ldots,a_i,\ldots,a_l']$, where $i<l\le n$ and $a_l'\le a_l$. By Corollary \ref{lem:9}, $[a_t,\ldots,a_i,\ldots,a_l']$ is a blowup of $[2,1,2]$ with exactly one entry equal to 1 and by Corollary \ref{4}, $[a_t,\ldots,a_i+1,\ldots,a_l']$ has dual that is a budding of $[4]$; denote this budding by $[b_1,\ldots,b_k]$. By Theorem \ref{contfrac}, $[a_t+1,\ldots,a_i+1,\ldots,a_l']$ has dual $[2,b_1,\ldots,b_k]$. By performing buddings to recover the original fraction, we have that the dual of $[a_1,\ldots,a_t+1,\ldots,a_i+1,\ldots,a_n]$ is a budding of $[2,b_1,\ldots,b_k]$, by Corollary \ref{cor:2}. If $t>i$, then a similar argument shows that the dual of $[a_1,\ldots,a_i+1,\ldots,a_t+1,\ldots,a_n]$ is a budding of $[b_1,\ldots,b_k,2]$, where $[b_1,\ldots,b_k]$ is a budding of $[4]$.

Now suppose $t=i$. The only such fraction of length 3 is $[2,3,2]$, which has dual $[3,3]$. Let $[a_1,\ldots,a_n]=0$ be a blowup of $[2,1,2]$ with exactly one entry equal to 1, namely $a_i$. By Corollary \ref{lem:9}, $[a_1,\ldots,a_n]$ is a budding of $[2,1,2]$. Thus $[a_1,\ldots,a_i+2,\ldots,a_n]$ is a budding of $[2,3,2]$. By Corollary \ref{cor:2}, the dual of $[a_1,\ldots,a_i+2,\ldots,a_n]$ must be a budding of the dual of $[2,3,2]$; that is, the dual is a budding of $[3,3]$.\end{proof}

\begin{lemma}[Type 2(ii)] Let $[a_1,\ldots,a_n]=0$ be a blowup of $[2,1,2]$ that is not a blowup of $[1,2,1]$ and suppose there are exactly two entries that are equal to 1, namely $a_i$ and $a_j$ (note, $i,j\notin\{1,n\}$ by Lemma \ref{lem:10}). Then $[a_1,\ldots,a_i+1,\ldots,a_j+1,\ldots,a_n]$ has dual that is a budding of $[2,b_1,\ldots,b_k,c_1,\ldots,c_l,2]$, where $[b_1,\ldots,b_k]$ and $[c_1,\ldots,c_l]$ are buddings of $[4]$. Moreover, this budding corresponds to a plumbing in Theorem \ref{thm:main}(d) \label{lemforthm4}\end{lemma}

\begin{proof} Note that the minimal length of such a continued fraction is 5 and the only such fraction of length 5 is $[3,1,3,1,3]=0$. Moreover, $[3,2,3,2,3]$ has dual $[2,4,4,2]$, which is the minimal length dual fraction of the desired form. We start by considering blowups of $[3,1,3,1,3]$. Notice that we can order the blowups so that we first perform all blowups before or after the first occurrence of 1 and then all blowups before or after the second occurrence of 1. By rewriting $[3,1,3,1,3]$ as $[1+\underline{2,1,2}+(-1)+\underline{2,1,2}+1]$, we can view the two underlined portions as subfractions $[2,1,2]$ of $[3,1,3,1,3]$. By repeatedly blowing up $[3,1,3,1,3]$ before or after the first occurrence of 1, we also blow up the first subfraction before or after 1. Thus we end up with a fraction of the form $[1+x_1,\ldots,x_i,\ldots, x_p+(-1)+2,1,3]$, where $[x_1,\ldots,x_p]=0$ is a blowup of $[2,1,2]$ with exactly one entry equal to 1, denoted by $x_i$. Now blow up before or after the second occurrence of 1 to obtain a fraction of the form $[1+x_1,\ldots,x_i,\ldots, x_p+(-1)+y_1,\ldots,y_j,\ldots,y_q+1]$, where $[y_1,\ldots,y_q]=0$ is a blowup of $[2,1,2]$ with exactly one entry equal to 1, denoted by $y_j$. Thus by Corollary \ref{4}, $[x_1,\ldots,x_i+1,\ldots,x_p]$ and $[y_1,\ldots,y_j+1,\ldots,y_q]$ have duals that are blowups of $[4]$, which we denote by $[b_1,\ldots,b_k]$ and $[c_1,\ldots,c_l]$, respectively. By Corollaries \ref{lem:13} and \ref{cor:2}, $[1+x_1,\ldots,x_i,\ldots, x_p+(-1)+y_1,\ldots,y_j,\ldots,y_q+1]$ has dual $[2,b_1,\ldots,b_k,c_1,\ldots,c_l,2]$. Thus any blowup of $[3,1,3,1,3]$ with exactly two entries equal to 1 has dual of the desired form.

Now let $[a_1,\ldots,a_n]$ be as in the statement of the lemma. Since $[a_1,\ldots,a_n]$ is not a blowup of $[1,2,1]$, the only way to obtain $[a_1,\ldots,a_n]$ is to first blow up $[2,1,2]$ before or after the middle entry. This yields $[3,1,2,2]$ or $[2,2,1,3]$. Furthermore, at each step, we cannot blow up at the beginning or end of the fraction. Otherwise, we would obtain a fraction with first or last entry equal to 1, which contradicts Lemma \ref{lem:10}. Thus either $a_1\ge 3$, $a_n\ge 3$, or $a_1,a_n\ge3$. If $a_1,a_n\ge3$, then we claim that $[a_1,\ldots,a_n]$ is a blowup of $[3,1,3,1,3]$. Notice that the only way to obtain such a fraction is by first blowing up a fraction of the form $[2^{[n-1]},1,n]$ (or $[n,1,2^{[n-1]}]$) between the first two entries (or last two entries) to obtain $[3,1,3,2^{[n-3]},1,n]$ (or $[n,1,2^{[n-3]},3,1,3]$). We can then blow up before or after either occurrence of 1. Thus such continued fractions are clearly blowups of $[3,1,3,1,3]$, which we already handled in the previous paragraph.

Now suppose $a_1=2$ and $a_n\ge 3$ (or similarly, $a_1\ge3$ and $a_n=2$) and let $a_i=a_j=1$, where $1<i,j<n$. We claim that by a sequence of debuddings, we can obtain a fraction that is a blowup of $[3,1,3,1,3]$. By Lemma \ref{lem:8}, the first debudding $[a_2,\ldots,a_n-1]$ is a blowup of $[2,1,2]$ and not a blowup of $[1,2,1]$. Moreover, it still has two entries equal to 1 and $a_2,a_n-1\ge 2$. If $a_2,a_n-1\ge3$, then by the previous paragraph, we are done. If $a_2=2$ and $a_n-1\ge3$ (or vice versa), then we can perform another debudding. Since the fraction has finite length, this process terminates, yielding a fraction with first and last entry greater than 2 (they cannot both be equal to 2 by the remarks above). Thus the result is a blowup $[3,1,3,1,3]$. Call this fraction $[a_1',\ldots,a_m']$, where $a_i'=a_j'=1$, $i\neq j$, and $i,j\neq1,n$. By the calculation in the first paragraph, $[a_1',\ldots,a_i'+1,\ldots,a_j'+1,\ldots,a_m']$ has dual of the form $[2,b_1,\ldots,b_k,c_1,\ldots,c_l,2]$, where $[b_1,\ldots,b_k]$ and $[c_1,\ldots,c_l]$ are buddings of $[4]$. Now we can perform buddings to $[a_1',\ldots,a_i'+1,\ldots,a_j'+1,\ldots,a_m']$ to recover the original fraction $[a_1,\ldots,a_i+1,\ldots,a_j+1,\ldots,a_n]$. Thus by Corollary \ref{cor:2}, its dual fraction can obtained by a sequence of buddings of $[2,b_1,\ldots,b_k,c_1,\ldots,c_l,2]$.\end{proof}

\bibliographystyle{plain}
\bibliography{Bibliography}

\end{document}